\documentclass[AMA,STIX1COL,article]{WileyNJD-v2}

\articletype{RESEARCH ARTICLE}%
\newtheorem{exam}{Example}
\newtheorem{lemm}{Lemma}
\usepackage{setspace}
\usepackage{orcidlink}
\usepackage{color}
\usepackage{graphicx}
\usepackage{float}
\usepackage{epstopdf}
\usepackage{epsfig}
\usepackage{subfigure}
\usepackage{hyperref}
\usepackage[T1]{fontenc}
\usepackage[titles,subfigure]{tocloft}
\hypersetup{hidelinks}
\received{}
\revised{}
\accepted{}
\begin{document}
\title{Dynamical behavior and optimal control of a stochastic SAIRS epidemic model with two saturated incidences}

\author[]{Xiaohui Zhang}

\author[]{Zhiming Li~\orcidlink{ https://orcid.org/0000-0002-3646-0784}}


\author[]{Shenglong Chen}
\author[]{Jikai Yang}

\authormark{Zhang \textsc{et al.}}

\address[]{ College of Mathematics and System Science, Xinjiang University,  Urumqi, China}

\corres{\email{zmli@xju.edu.cn.}}

\presentaddress{Natural Science Foundation, China, Grant Number: 12061070,
Xinjiang Uygur Autonomous Region Natural Science Foundation, China, Grant Number: 2021D01E13, Xinjiang Uygur Autonomous Region Postgraduates Research Innovation Program Foundation, China, Grant Number: XJ2023G016, XJ2023G017,
2023 Annual Planning Project of Commerce Statistical Society of China, Grant Number: 2023STY61,
Innovation Project of Excellent Doctoral Students of Xinjiang University, China, Grant Number:  XJU2023BS017.}

\abstract{Stochastic models are widely used to investigate the spread of epidemics in a complex environment. This paper extends a deterministic SAIRS epidemic model  \textcolor{blue}{(Ottaviano et al. 2022)}  to a stochastic case with limited patient capacity and exposure.  We first study the dynamical properties of the model under certain conditions, including persistence, extinction, and ergodic. Then, we introduce vaccination and isolation into the model as control variables. The optimal control strategies are obtained based on Pontryagin's minimum principle. Finally, numerical simulations are given to illustrate our theoretical results.}

\maketitle

\section{Introduction}


Infectious diseases are one of the significant factors which threaten people's health. How to effectively prevent and control the spread of an epidemic has been one of the hottest topics worldwide. Mathematical modeling often plays a vital role in the dynamics of various infectious diseases, such as the deterministic SIR model and extended versions \textcolor{blue}{(Liu and Chen 2016; Gao and Zhang et al. 2019; Acemoglu et al. 2021; Basnarkov et al. 2022; Li and Eskandari 2023; Khan et al. 2024; Khan and Khan 2024).}  Due to environmental noise or other factors, many scholars have made an in-depth study and research of stochastic models \textcolor{blue}{ (Gray et al. 2011; Zhang and Wang 2014; Cai et al. 2015; Lei and Yang 2017; Liu et al. 2017;  Tain et al. 2020;   Din et al. 2022;  Din and Li 2022;  Wang et al. 2023).}  Recently,  \textcolor{blue}{ Ottaviano et al. (2022)} proposed a deterministic SAIRS model with asymptomatic individuals. Suppose that the total population $N$ is divided into four compartments: susceptible (S), asymptomatic (A), infected (I), and recovered (R) individuals. Let $S(t), A(t), I(t)$, and $R(t)$ be the numbers of the susceptible, asymptomatic, infected, and recovered individuals, respectively. The SAIRS model with two bilinear incidences is defined by
\begin{equation}\label{eq0}
\left\{\begin{aligned}
d S(t)& =[\Lambda-(\beta_{A} A(t)+\beta_{I} I(t)) S(t)-(\mu+\nu) S(t)+\gamma R(t)]dt, \\
d A(t)& =[(\beta_{A} A(t)+\beta_{I} I(t))S(t)-\left(\alpha+\delta_{A}+\mu\right) A(t)]dt, \\
d I(t)& =[\alpha A(t)-\left(\delta_{I}+\mu\right) I(t)]dt, \\
d R(t) & =[\delta_{A} A(t)+\delta_{I} I(t)+\nu S(t)-(\gamma+\mu) R(t)]dt,
\end{aligned}
\right.
\end{equation}
where $\Lambda$ is the migration rate of the susceptible individuals, $\beta_{A}$ and $\beta_{I}$ are the infection rates of asymptomatic and symptomatic patients, $\mu$ is the natural birth rate,  $\nu$ is the rate of acquisition of immunity after vaccination, $\delta_{A} $ and  $\delta_{I} $ represent natural recovery rates of asymptomatic and infected individuals, $\alpha$ is the transformation rate from asymptomatic to patient, and $\gamma $ is the transformation rate from recovered to susceptible. Here, all parameters are positive in the model (\ref{eq0}).

For the SAIRS model (\ref{eq0}), one of the advantages is based on the fact that both asymptomatic and infectious individuals can infect susceptible ones. However, there are also some limitations: (i) It ignores the effect of an upper bound on the number of person-to-person contacts during epidemic transmission. In the model (\ref{eq0}), the bilinear incidence relatively depends on the number of susceptible individuals to the infection. This is because a higher number of susceptible individuals increases the likelihood of asymptomatic or severely infected people coming into contact with them, leading to an increased incidence \textcolor{blue}{ (Wang  et al. 2010;  Feng et al. 2022)}. When the number of susceptible individuals to the infection is substantial, a patient's capacity to interact with others is inevitably limited  \textcolor{blue}{ (Idoon et al. 2022)}. The incidence will reach the upper bound because of the limited exposure of those with the disease to the susceptible population. Thus, it is not reasonable that the incidence is proportional to the number of susceptible people. The saturated incidence rate considers the infected population's behavioral changes and crowding effects. It prevents the unboundedness of the exposure rate by choosing appropriate parameters, thus making the model (\ref{eq0}) more reasonable \textcolor{blue}{(Loyinmi et al. 2021).}  (ii) The influence of the diseased death rate should not be neglected. Without the diseased death rate, the model (\ref{eq0}) is only feasible to study the spread of an epidemic in a short time, but not a long time \textcolor{blue}{(Sempe et al. 2021)}.
Based on this point,  a modified SAIRS model with two saturated incidences is established as follows
\begin{equation}\label{eq1}
\left\{\begin{array}{l}
d S(t)=\left[\Lambda -\left(\frac{\beta_{I} I(t)}{1+b I(t)}+\frac{\beta_{A} A(t)}{1+b A(t)}\right) S(t)-\mu S(t)+\gamma R(t)\right]dt,\\
d A(t)=\left[\left(\frac{\beta_{I} I(t)}{1+b I(t)}+\frac{\beta_{A} A(t)}{1+b A(t)}\right) S(t)-\left(\alpha+\delta_{A}+\mu\right) A(t)\right]dt, \\
d I(t)=[\alpha A(t)-\left(\delta_{I}+\mu+d\right) I(t)]dt, \\
d R(t)=[\delta_{A} A(t)+\delta_{I} I(t)-(\gamma+\mu) R(t)]dt,
\end{array}
\right.
\end{equation}
where $d(\geq0) $ is the diseased death rate, and $b(\geq 0)$ is a coefficient of saturation incidence.  If $d=b=0$, then the model (\ref{eq1}) is equivalent to the model (\ref{eq0}) without vaccination.

In practice, there are various stochastic environmental factors, such as climate and politics, during the spread of an epidemic.   Through the analysis above,  it is necessary to consider  the stochastic version of the deterministic model (\ref{eq1}) as follows
\begin{equation}\label{eq2}
\left\{\begin{array}{l}
		d S(t)=\left[\Lambda -\left(\frac{\beta_{I} I(t)}{1+b I(t)}+\frac{\beta_{A} A(t)}{1+b A(t)}\right) S(t)-\mu S(t)+\gamma R(t)\right]dt+\sigma_{1} S(t) d B_{1}(t),\\
		d A(t)=\left[\left(\frac{\beta_{I} I(t)}{1+b I(t)}+\frac{\beta_{A} A(t)}{1+b A(t)}\right)S(t)-\left(\alpha+\delta_{A}+\mu\right) A(t)\right]dt+\sigma_{2}A(t) d B_{2}(t), \\
		d I(t)=[\alpha A(t)-\left(\delta_{I}+\mu+d\right) I(t) ]dt+\sigma_{3} I(t) dB_{3}(t),\\
		d R(t)=[\delta_{A} A(t)+\delta_{I} I(t)-(\gamma+\mu) R(t)]dt+\sigma_{4}R(t) d B_{4}(t),
\end{array}
\right.
\end{equation}
where  $B_{i}(t) (i=1,2,3,4)$  are independent Brownian motion and  $\sigma_{i} (i=1,2,3,4)$  denote the intensity of the white noise. Although vaccination has been considered in the model (\ref{eq0}), the mutations of some infectious diseases often make it invalid. Given the immense power of modern biotechnology, new vaccines become possible as an effective control strategy \textcolor{blue}{(Makinde et al. 2007; Din and Li 2021;  Din. 2021;  Din and Qura 2022; Dwivedi et al. 2022; Mondal et al. 2022; Bentaleb et al. 2023; Silver et al. 2023;  Zhang et al. 2024)}.  Based on this,  two control variables will be introduced in the stochastic model, such as vaccination and isolation of infected and asymptomatic patients, into the model (\ref{eq2}). Therefore,  the novelties of this paper are outlined as follows: (i) the stochastic SAIR model (\ref{eq2}) is more suitable than the deterministic model (\ref{eq1}) in a complex environment. (ii) Unlike the bilinear incidence, saturated incidences consider the effect of the population, not just the susceptible individuals. (iii)  Optimal control strategies are given to prevent the spread of the epidemic.

The rest of this paper is organized as follows.  In Section 2,  we review the basic concepts and useful lemmas needed in work.  In Section 3, we study the persistence, extinction, and stationary distribution of the stochastic SAIRS model (\ref{eq2}). In Section 4,  we establish a stochastic control system and obtain optimal strategies. To illustrate our theoretical results, numerical simulation is provided in Section 5.  Finally, we provide a summary in Section 6.

\section{Preliminaries}
Let $X(t)=(X_1(t),\ldots,X_n(t))$ be an $n$-dimensional stochastic process,  satisfing the following stochastic differential equation
\begin{eqnarray}\label{sde}
d X(t)=f(X(t), t) d t+g(X(t), t) d B(t), \quad t \geq t_{0},
\end{eqnarray}
with initial value  $X(0)=X_{0} \in \mathbb{R}^{n}$, where  $B(t)$  represents an $n$-dimensional standard Brownian motion. Denote by  $C^{2,1}\left(\mathbb{R}^{n} \times\left[t_{0},+\infty\right) ; \mathbb{R}_{+}\right)$ the family of all nonnegative functions $ V(X, t)$  defined on  $\mathbb{R}^{n}\times\left[t_{0}+\infty\right)$  such that they are twice continuously differentiable for $X$,  and once in  $t$. The operator $L$ for the system (\ref{sde}) is defined as follows
$$
L=\frac{\partial}{\partial t}+\sum_{i=1}^{n} f_{i}(X, t) \frac{\partial}{\partial X_{i}}+\frac{1}{2} \sum_{i, j=1}^{n}\left[g^{T}(X, t) g(X, t)\right]_{i j} \frac{\partial^{2}}{\partial X_{i} \partial X_{j}}.
$$
If $L$ acts on the function $V \in C^{2,1}\left(\mathbb{R}^{n} \times\left[t_{0},+\infty\right) ; \mathbb{R}_{+}\right)$, then
$$
L V(X, t)=V_{t}(X, t)+V_{X}(X, t) f(X, t)+\frac{1}{2} \operatorname{trace}\left[g^{T}(X, t) V_{X X}(X, t) g(X, t)\right],
$$
where  $$V_{t}=\frac{\partial V}{\partial t}, V_{X}=\left(\frac{\partial V}{\partial X_{1}}, \frac{\partial V}{\partial X_{2}}, \ldots, \frac{\partial V}{\partial X_{n}}\right), V_{X X}=\left(\frac{\partial^{2} V}{\partial X_{i} \partial X_{j}}\right)_{n \times n}.$$
 According to Ito's formula,
$$
d V(X(t), t)=L V(X(t), t) d t+V_{X}(X(t), t) g(X(t), t) d B(t), \quad X(t) \in \mathbb{R}^{n} .
$$

For convenience, this article presents the following two theorems that will be used in later sections.
\begin{lemm}\label{L1}
 Let $ (S(t), A(t), I(t), R(t))$  be the solution of model (\ref{eq2}) with initial value $(S(0), A(0), I(0), R(0)) \in R_{+}^{4}$ . Then,
$$
\limsup _{t \rightarrow+\infty}(S(t)+A(t)+I(t)+R(t))<+\infty \text { a.s. }
$$
and
\begin{equation}
\lim _{t \rightarrow+\infty} \frac{\int_{0}^{t} S(r) d B_{1}(r)}{t}= \lim _{t \rightarrow+\infty} \frac{\int_{0}^{t} A(r) d B_{2}(r)}{t}=
\lim _{t \rightarrow \infty} \frac{\int_{0}^{t} I(r) d B_{3}(r)}{t}= \lim _{t \rightarrow+\infty} \frac{\int_{0}^{t} R(r) d B_{4}(r)}{t}=0, \quad \text { a.s. }
\end{equation}
\end{lemm}
The proof is similar to that of \textcolor{blue}{ Taki et al. (2022)}, so it's omitted here.
\begin{lemm}
\label{L2} (Khasminskii Theorem \textcolor{blue}{(Khasminskii 2012)}). The Markov process  $X(t)$  has a unique stationary distribution  $\pi(\cdot)$  if a bounded domain  $U_{\varepsilon} \subset \mathbb{R}_{+}^{4}$  with regular boundary  $\Gamma$  exists and

(i) there exists a constant  $M>0$  satisfying    $\sum_{i, j=1}^{d} a_{i j}(x) \xi_{i} \xi_{j} \geq M|\xi|^{2}$, $X(0)=x \in U_{\varepsilon}$, $\xi \in \mathbb{R}_{+}^{4}$ .

(ii) there is a  $C^{2}$ -function  $V \geq 0 $ such that   $L V $  is negative for any $ \mathbb{R}_{4}^{+} \backslash U_{\varepsilon}$ . Then
$$
\mathbb{P}\left\{\lim _{T \rightarrow+\infty} \frac{1}{T} \int_{0}^{T} f(X(t)) \mathrm{d} t=\int_{\mathbb{R}_{+}^{4}} f(x) \pi(\mathrm{d} x)\right\}=1
$$
for all  $X\in \mathbb{R}_{+}^{4}$ , where  $f(\cdot)$  is an integrable function with respect to the measure  $\pi$ .
\end{lemm}
\section{ Dynamical properties}
 In this section,  this article mainly studies the dynamic properties of the model (\ref{eq2}), including the solution's persistence, extinction, and stationary distribution.
\subsection{ Stochastic persistence}

In this subsection, we focus on the persistence of the solution in the model  (\ref{eq2}).  Let $f(t)$ be a continuous integrable function. The function $f(t)$ is called stochastic persistence if  $\lim _{t \rightarrow \infty} \frac{1}{t} \int_{0}^{t} f(s) \mathrm{d} s>0$. Denote $m=4\mu+\alpha+\delta_{A}+\delta_{I}+d+\gamma+\frac{1}{2}\sum_{i=1}^4\sigma_i^2,$
and
$$
R_{0}^{s}:=\frac{3}{m}\sqrt[3]{\Lambda\alpha\beta}.$$

\begin{theorem}\label{th2.2}
 If $R_{0}^{s}>1$,  for any initial value  $(S(0), A(0), I(0), R(0))\in \mathbb{R}_{+}^{4}$, the solution of system (\ref{eq2}) satisfies the following inequalities
\begin{equation*}
\begin{aligned}
&\lim _{t \rightarrow \infty} \frac{1}{t} \int_0^t S(s) \mathrm{d} s>\frac{mb(\mu+\alpha+\delta_{A})(\mu+d+\delta_{I})(R_{0}^{s}-1)}{(\beta_{A}+\beta_{I})(\alpha\beta_{I}+(\beta_{A}+b\alpha)(\delta_{I}+\mu+d))}>0,\\
&\lim _{t \rightarrow \infty} \frac{1}{t} \int_{0}^{t} A(s) \mathrm{d} s>\frac{m(\delta_{I}+\mu+d)(R_{0}^{s}-1)}{\alpha\beta_{I}+(\beta_{A}+b\alpha)(\delta_{I}+\mu+d)}>0,\\
&\lim _{t \rightarrow \infty} \frac{1}{t} \int_{0}^{t} I(s) \mathrm{d} s > \frac{{m\alpha(R_{0}^{s}-1)}}{\alpha\beta_{I}+(\beta_{A}+b\alpha)(\delta_{I}+\mu+d)}>0,\\
&\lim _{t \rightarrow \infty} \frac{1}{t} \int_{0}^{t} R(s) \mathrm{d} s>\frac{m(\delta_{A}(\mu+d+\delta_{I})+\alpha\delta_{I})(R_{0}^{s}-1)}{(\gamma+\mu)(\alpha\beta_{I}+(\beta_{A}+b\alpha)(\delta_{I}+\mu+d))}>0.
\end{aligned}
\end{equation*}
That is, the solution of the system is persistent.
\end{theorem}
\begin{proof}
Construct a function $V(S,A,I,R) = -\ln S - \ln A - \ln I - \ln R$. Using Ito's formula \textcolor{blue}{(Mao 2007)}, it follows that
\begin{equation}\label{15}\nonumber
\mathrm{d} V_{1}(t)=L V_{1}(t) \mathrm{d} t-\sigma_{1} \mathrm{~d} B_{1}(t)-\sigma_{2} \mathrm{~d} B_{2}(t)-\sigma_{3} \mathrm{~d} B_{3}(t)-\sigma_{4} \mathrm{~d} B_{4}(t),
\end{equation}
where
\begin{equation}\nonumber
\begin{aligned}
L V_{1}(t) =&-\frac{1}{S(t)}\left(\Lambda -(\frac{\beta_{I} I(t)}{1+b I(t)}+\frac{\beta_{A} A(t)}{1+b A(t)}) S(t)-\mu S(t)+\gamma R(t)\right)\\
&-\frac{1}{A(t)}((\frac{\beta_{I} I(t)}{1+b I(t)}+\frac{\beta_{A} A(t)}{1+b A(t)})S(t)-\left(\alpha+\delta_{A}+\mu\right) A(t))\\
&-\frac{1}{I(t)}\left(\alpha A(t)-\left(\delta_{I}+\mu+d\right) I(t)\right)+\frac{1}{2}\left(\sigma_{1}^{2}+\sigma_{2}^{2}+\sigma_{3}^{2}+\sigma_{4}^{2}\right)\\
&-\frac{1}{R(t)}\left(\delta_{A} A(t)+\delta_{I} I(t)-(\gamma+\mu) R(t)\right) \\
=&-\frac{\Lambda}{S(t)}-\gamma\frac{R(t)}{S(t)}-\frac{\beta_{I} I(t)S(t)}{(1+b I(t))A(t)}-\frac{\beta_{A} S(t) }{(1+b A(t))}-\alpha\frac{A(t)}{I(t)}-\delta_{A}\frac{A(t)}{R(t)}-\delta_{I}\frac{I(t)}{R(t)}\\
&+m+\frac{\beta_{I} I(t)}{1+b I(t)}+\frac{\beta_{A} A(t)}{1+b A(t)}\\
\leq&-\frac{\Lambda}{S(t)}-\beta(\frac{ I(t)S(t)}{(1+b I(t))A(t)}+\frac{ S(t) }{(1+b A(t))})-\alpha\frac{A(t)}{I(t)}-\gamma\frac{R(t)}{S(t)}-\delta_{A}\frac{A(t)}{R(t)}-\delta_{I}\frac{I(t)}{R(t)}\\
&+m+\frac{\beta_{I} I(t)}{1+b I(t)}+\frac{\beta_{A} A(t)}{1+b A(t)}\\
<&-\frac{\Lambda}{S(t)}-\frac{\beta I(t)S(t)}{(1+b I(t))A(t)}-\alpha\frac{A(t)(1+b I(t))}{I(t)}+m+b\alpha A(t)+\frac{\beta_{I} I(t)}{1+b I(t)}+\frac{\beta_{A} A(t)}{1+b A(t)},\\
\end{aligned}
\end{equation}
where $\beta=\min\left\{\beta_{I},\beta_{A}\right\}$.
 Because of the inequality $a+b+c\ge3\sqrt[3]{abc}$, combining the above inequalities yields
\begin{equation}\label{16}
\begin{aligned}
L V_{1}(t) & <-3\sqrt[3]{\Lambda\alpha\beta}+\beta_{I} I(t)+\beta_{A} A(t)+m=-m(R_{0}^{s}-1)+\beta_{I} I(t)+(\beta_{A}+b\alpha) A(t),
\end{aligned}
\end{equation}
 Integrate both sides of the above equation (\ref{16}) and divide by $t$. As $t$ tend to infinity, it leads that
\begin{eqnarray}\label{17}
0 &<&- m(R_{0}^{s}-1)+\lim _{t \rightarrow \infty} \frac{\beta_{I}}{t} \int_{0}^{t} I(s) \mathrm{d} s+\lim _{t \rightarrow \infty} \frac{\beta_{A}+b\alpha}{t} \int_{0}^{t} A(s) \mathrm{d} s \nonumber\\
&& -\lim _{t \rightarrow \infty} \frac{1}{t} \int_{0}^{t} \sigma_{1} \mathrm{~d} B_{1}(s)-\lim _{t \rightarrow \infty} \frac{1}{t} \int_{0}^{t} \sigma_{2} \mathrm{~d} B_{2}(s)-\lim _{t \rightarrow \infty} \frac{1}{t} \int_{0}^{t} \sigma_{3} \mathrm{~d} B_{3}(s)-\lim _{t \rightarrow \infty} \frac{1}{t} \int_{0}^{t} \sigma_{4} \mathrm{~d} B_{4}(s).
\end{eqnarray}
According to Lemma \ref{L1}, we can derive
\begin{equation}\nonumber
\lim _{t \rightarrow \infty} \frac{1}{t} \int_{0}^{t} \sigma_{1} \mathrm{~d} B_{1}(s)=\lim _{t \rightarrow \infty} \frac{1}{t} \int_{0}^{t} \sigma_{2} \mathrm{~d} B_{2}(s)=\lim _{t \rightarrow \infty} \frac{1}{t} \int_{0}^{t} \sigma_{3} \mathrm{~d} B_{3}(s)=\lim _{t \rightarrow \infty} \frac{1}{t} \int_{0}^{t} \sigma_{4} \mathrm{~d} B_{4}(s)=0.
\end{equation}
So (\ref{17}) becomes
\begin{equation}\label{18}
\lim _{t \rightarrow \infty} \frac{\beta_{I}}{t} \int_{0}^{t} I(s) \mathrm{d} s+\lim _{t \rightarrow \infty} \frac{\beta_{A}+\beta\alpha}{t} \int_{0}^{t} A(s) \mathrm{d} s \geq m(R_{0}^{s}-1).
\end{equation}
By the third equation of the system (\ref{eq2}), integrating both sides of the equation. Thus,
\begin{equation}\nonumber
I(t)-I(0)=\int_{0}^{t}[\alpha A(s)-(\mu+d+\delta_{I}) I(s)] \mathrm{d} s+\int_{0}^{t} \sigma_{3} I(s) \mathrm{d} B_{3}(s).
\end{equation}
Both sides are simultaneously divided by  $t$ and take the limit.  It is obvious that
\begin{equation}\label{19}
\lim _{t \rightarrow \infty} \frac{\mu+d+\delta_{I}}{t} \int_{0}^{t} I(s) \mathrm{d} s=\lim _{t \rightarrow \infty} \frac{\alpha}{t} \int_{0}^{t} A(s) \mathrm{d} s+\lim _{t \rightarrow \infty} \frac{1}{t} \int_{0}^{t} \sigma_{3} I(s) \mathrm{d} B(s).
\end{equation}
Based on the strong law of large numbers for martingales\cite{bib7},  the equation (\ref{19}) becomes
\begin{equation}\nonumber
\lim _{t \rightarrow \infty} \frac{\delta_{I}+d+\mu}{t} \int_{0}^{t} I(s) \mathrm{d} s=\lim _{t \rightarrow \infty} \frac{\alpha}{t} \int_{0}^{t} A(s) \mathrm{d} s.
\end{equation}
Thus,
\begin{equation}\label{20}
\lim _{t \rightarrow \infty} \frac{1}{t} \int_{0}^{t} A(s) \mathrm{d} s=\lim _{t \rightarrow \infty} \frac{\delta_{I}+\mu+d}{\alpha} \cdot \frac{1}{t} \int_{0}^{t} I(s) \mathrm{d} s.
\end{equation}
According to  (\ref{20}) and (\ref{18}), the following inequality hold:
\begin{equation}\nonumber
\lim _{t \rightarrow \infty} \frac{1}{t} \int_{0}^{t} I(s) \mathrm{d} s > \frac{{\alpha m(R_{0}^{s}-1)}}{\alpha\beta_{I}+(\beta_{A}+b\alpha)(\delta_{I}+\mu+d)},
\\
\lim _{t \rightarrow \infty} \frac{1}{t} \int_{0}^{t} A(s) \mathrm{d} s>\frac{(\delta_{I}+\mu+d)m(R_{0}^{s}-1)}{\alpha\beta_{I}+(\beta_{A}+b\alpha)(\delta_{I}+\mu+d)}.
\end{equation}
By the second equation of the system (\ref{eq2}),  it follows that
\begin{equation*}\label{40}
\begin{aligned}
&\mathrm{d} A \leq (\frac{\beta_{A}+\beta_{I}}{b}S(t)-(\mu+\alpha+\delta_{A})A(t))dt+\sigma_{2} A(t) \mathrm{~d} B_{2}(t),\\
&\mathrm{d} R = \delta_{A}A(t)+\delta_{I}I(t)-(\mu+\gamma)R(t))dt+\sigma_{4} R(t) \mathrm{~d} B_{4}(t).
\end{aligned}
\end{equation*}
Integrate the above equations from 0 to $t$ and divide $t$ on both sides. Thus, the following results are obtained
\begin{equation*}\label{22}
\begin{aligned}
&\lim _{t \rightarrow \infty} \frac{1}{t} \int_{0}^{t} S(s) \mathrm{d} s>\frac{b(\mu+\alpha+\delta_{A})(\mu+d+\delta_{I})m(R_{0}^{s}-1)}{(\beta_{A}+\beta_{I})(\alpha\beta_{I}+(\beta_{A}+b\alpha)(\delta_{I}+\mu+d))},\\
&\lim _{t \rightarrow \infty} \frac{1}{t} \int_{0}^{t} R(s) \mathrm{d} s>\frac{(\delta_{A}(\mu+d+\delta_{I})+\alpha\delta_{I})m(R_{0}^{s}-1)}{(\gamma+\mu)(\alpha\beta_{I}+(\beta_{A}+b\alpha)(\delta_{I}+\mu+d))}.
\end{aligned}
\end{equation*}The proof is complete.
\end{proof}

\subsection{Extinction of infectious disease}
\begin{theorem}\label{th2.3}
If $\beta=\max\left\{\beta_{I},\beta_{A}\right\}$, $ h=\min \left\{\mu+\delta_{A}+\frac{\sigma_{2}^{2}}{2},(\mu+\delta_{I}+d)+\frac{\sigma_{3}^{2}}{2}\right\}$, and $\beta \frac{\Lambda}{\mu}-\frac{h}{2}<0$, for any initial value $(S(0), A(0), I(0), R(0))\in \mathbb{R}_{+}^{4}$,  then
\begin{equation}\nonumber
\lim_{t\to\infty} S(t) = \frac{\Lambda}{\mu},\quad
\lim _{t \rightarrow \infty} A(t=\lim _{t \rightarrow \infty} I(t)=\lim _{t \rightarrow \infty} R(t)=0.
\end{equation}
\end{theorem}
\begin{proof}
 Define a function  $V_{3}=\ln (A(t)+I(t))$. According to Ito's formula, the
\begin{equation}\label{23}
\mathrm{d} V_{3}(t)=L V_{3}(t) \mathrm{d} t+\frac{\sigma_{2} A}{A+I} \mathrm{~d} B_{2}(t)+\frac{\sigma_{3} I}{A+I} \mathrm{~d} B_{3}(t),
\end{equation}
where
\begin{eqnarray}\label{24}
L V_{3} & =&\frac{1}{A+I}\left[\frac{\beta_{I} S I}{1+bI}+\frac{\beta_{A} S A}{1+bA}-(\mu+\delta_{A}+\alpha) A+\alpha A-(\mu+\delta_{I}+d) I\right]-\frac{\sigma_{2}^{2} A^{2}+\sigma_{3}^{2} I^{2}}{2(A+I)^{2}} \nonumber\\
& =&\frac{\frac{\beta_{I} S I}{1+bI}+\frac{\beta_{A} S A}{1+bA}}{A+I}+\frac{\left(-\mu-\delta_{A}-\frac{\sigma_{2}^{2}}{2}\right) A^{2}+\left[-(\mu+\delta_{I}+d)-\frac{\sigma_{3}^{2}}{2}\right] I^{2}-(2\mu+d+\delta_{A}+\delta_{I}) A I}{(A+I)^{2}} \nonumber\\
& \leq& \beta S+\frac{\left(-\mu-\delta_{A}-\frac{\sigma_{2}^{2}}{2}\right) A^{2}+\left[-(\mu+\delta_{I}+d)-\frac{\sigma_{3}^{2}}{2}\right] I^{2}}{(A+I)^{2}} \nonumber\\
& \leq& \beta S-\frac{\left(\mu+\delta_{A}+\frac{\sigma_{2}^{2}}{2}\right) A^{2}+\left[(\mu+\delta_{I}+d)+\frac{\sigma_{3}^{2}}{2}\right] I^{2}}{2\left(A^{2}+I^{2}\right)}\nonumber \\
& \leq& \beta S-h/2.
\end{eqnarray}
 Combine (\ref{24}) and (\ref{23}),  integrate from 0 to $t$ and divide $t$ on both sides, then take the limit. Thus,
\begin{equation}\label{25}
\lim _{t \rightarrow \infty} \frac{\ln (A+I)}{t} \leq \lim _{t \rightarrow \infty}\frac{\beta}{t} \int_{0}^{t} S(s) \mathrm{d} s-\frac{h}{2}.
\end{equation}
On the other hand, the four equations for the system (\ref{eq2}) are added together. Hence,
\begin{equation}\label{41}
\mathrm{d}(S+A+I+R)=[\Lambda-\mu(S+A+I+R)-d I] \mathrm{d} t+\sigma_{1} S \mathrm{~d} B_{1}+\sigma_{2} A \mathrm{~d} B_{2}+\sigma_{3} I \mathrm{~d} B_{3}+\sigma_{4} I \mathrm{~d} B_{4}.
\end{equation}
Integrate the equation (\ref{41}) and divide by $t$.  By applying the strong law of large numbers,  it is obvious that
\begin{eqnarray*}
\lim _{t \rightarrow \infty} \frac{1}{t} \int_{0}^{t} \sigma_{1} S(s) \mathrm{d} B(s)=\lim _{t \rightarrow \infty} \frac{1}{t} \int_{0}^{t} \sigma_{2} A(s) \mathrm{d} B(s)=\lim _{t \rightarrow \infty} \frac{1}{t} \int_{0}^{t} \sigma_{3} I(s) \mathrm{d} B(s)=\lim _{t \rightarrow \infty} \frac{1}{t} \int_{0}^{t} \sigma_{4} R(s) \mathrm{d} B(s)=0.
\end{eqnarray*}
Hence,
\begin{equation}\nonumber
0 \leq \Lambda-\lim _{t \rightarrow \infty} \frac{\mu}{t} \int_{0}^{t}(S(s)+A(s)+I(s)+R(s))-dI(s) \mathrm{d} s.
\end{equation}
That is to say,
\begin{equation}\label{26}
\lim _{t \rightarrow \infty} \frac{1}{t} \int_{0}^{t} S(s) \mathrm{d} s \leq\lim _{t \rightarrow \infty} \frac{1}{t} \int_{0}^{t}(S(s)+A(s)+I(s)+R(s))-dI(s) \mathrm{d} s\leq\frac{\Lambda}{\mu}.
\end{equation}
Combining (\ref{26}) and (\ref{25}) leads to
\begin{equation}\nonumber
\lim _{t \rightarrow \infty} \frac{\ln (A+I)}{t} \leq\frac{\beta \Lambda}{\mu}-\frac{h}{2}<0.
\end{equation}
So,
$
\lim _{t \rightarrow \infty} A(s)=\lim _{t \rightarrow \infty} I(s)=0.
$
Through the fourth equation in the system (\ref{eq2}),  $\lim _{t \rightarrow \infty} R(s)=0$.
Moreover, according to system (\ref{eq2}), for any small $\varepsilon > 0$, there exists a positive constant $T (> 0)$ such that $\frac{\beta_{I} I(t)}{1+b I(t)} + \frac{\beta_{A} A(t)}{1+b A(t)}<\varepsilon$ for $t > T$. Thus,
\begin{equation*}\label{4}
  \begin{aligned}
d S(t) \geq & \left(\Lambda-\left(\frac{\beta_{I} I(t)}{1+b I(t)}+\frac{\beta_{A} A(t)}{1+b A(t)}\right)S(t)-\mu S(t)\right) d t+\sigma_{1} S(t) d B_{1}(t) \\
\geq& \left(\Lambda-\mu S(t)-\varepsilon S(t)\right) d t+\sigma_{1} S(t) d B_{1}(t).
\end{aligned}
\end{equation*}
Integrating both sides of the above equation  and dividing by $t$,  one has
$
\liminf_{t\to\infty}S(t) \geq \frac{\Lambda}{\mu+\varepsilon}.
$
For any   sufficiently small and positive number $\varepsilon$,
$\liminf_{t\to\infty} S(t) \geq \Lambda/\mu$.
Integrate both sides of (\ref{41}) from 0 to $t$.  By the strong law of large numbers,  it follows that
\begin{equation*}\label{6}
  \limsup_{t\to \infty}(S(t)+A(t)+I(t)+R(t)-dI(t))\leq \Lambda/\mu.
\end{equation*}
It can be deduced that $\lim_{t\rightarrow+\infty} S(t) =\Lambda/\mu.$ The proof is complete.
\end{proof}

\subsection{ Stationary distribution and ergodic}
The stationary distribution plays a crucial role in research on stochastic epidemic models.
 By analyzing the ergodic property, the characteristics of the epidemic can be more effectively comprehended. Denote $a\vee b=\max\{a,b\}$.
\begin{theorem}\label{th2.4}
 If  $R_{0}^{s}>1$, for any initial value  $(S(0), A(0), I(0), R(0))\in \mathbb{ R}_{+}^{4}$, then system (\ref{eq2}) has a unique stationary distribution, and this distribution is ergodic.
\end{theorem}

\begin{proof}
Define a function
\begin{eqnarray*}
\widetilde{W}=W_{1}+W_{2}+W_{3}+W_{4},
\end{eqnarray*}
where  $W_{1}=\frac{1}{1+\theta}(S+A+I+R)^{1+\theta}, W_{2}=-\ln S-\ln A-\ln I,  W_{3}=-\lambda I,  W_{4}=-\ln S-\ln I-\ln R$. Here, $\lambda $ is a pending positive number, and $0<\theta<\frac{2\mu}{\sigma_{1}^{2}\vee \sigma_{2}^{2} \vee \sigma_{3}^{2}\vee \sigma_{4}^{2}}$.
Through calculation, $\widetilde{W}$ has a minimum value of $W_{0}$.
Define a Lyapunov function $ W=\widetilde{W}-W_{0}$. Denote  $H=\sup _{t \in(0, \infty)} \Lambda(S+A+I+R)^{\theta}-\frac{B}{2}(S+A+I+R)^{1+\theta}, B=\mu-\frac{\theta}{2}(\sigma_{1}^{2} \vee \sigma_{2}^{2} \vee \sigma_{3}^{2} \vee\sigma_{4}^{2}).$    By Ito's formula,
\begin{equation*}\label{27}
\begin{aligned}
L W_{1}=&  (S+A+I+R)^{\theta}[\Lambda-\mu(S+A+I+R)-dI]+\frac{\theta}{2}(S+A+I+R)^{\theta-1} \sigma_{1}^{2} S^{2} \\
& +\frac{\theta}{2}(S+A+I+R)^{\theta-1} \sigma_{2}^{2} A^{2}+\frac{\theta}{2}(S+A+I+R)^{\theta-1} \sigma_{3}^{2} I^{2} \\
&+\frac{\theta}{2}(S+A+I+R)^{\theta-1} \sigma_{3}^{2} R^{2}\\
\leq & \Lambda(S+A+I+R)^{\theta}-\mu(S+A+I+R)^{1+\theta}+\frac{\theta}{2}(S+A+I+R)^{1+\theta}\left(\sigma_{1}^{2} \vee \sigma_{2}^{2} \vee \sigma_{3}^{2}\vee \sigma_{4}^{2}\right) \\
= & \Lambda(S+A+I+R)^{\theta}-\left(\mu-\frac{\theta}{2}\left(\sigma_{1}^{2} \vee \sigma_{2}^{2} \vee \sigma_{3}^{2}\vee \sigma_{4}^{2}\right)\right)(S+A+I+R)^{1+\theta} \\
\leq & H-\frac{B}{2}(S^{1+\theta}+A^{1+\theta}+I^{1+\theta}+R^{1+\theta}),\\
L W_{2}(t) =&\frac{1}{S(t)}\left(\Lambda -(\frac{\beta_{I} I(t)}{1+b I(t)}+\frac{\beta_{A} A(t)}{1+b A(t)}) S(t)-\mu S(t)+\gamma R(t)\right)\\
&-\frac{1}{A(t)}((\frac{\beta_{I} I(t)}{1+b I(t)}+\frac{\beta_{A} A(t)}{1+b A(t)})S(t)-\left(\alpha+\delta_{A}+\mu\right) A(t)) I(t))\\
&-\frac{1}{I(t)}\left(\alpha A(t)-\left(\delta_{I}+\mu+d\right) I(t)\right)+\frac{1}{2}\left(\sigma_{1}^{2}+\sigma_{2}^{2}+\sigma_{3}^{2}+\sigma_{4}^{2}\right)\\
&-\frac{1}{R(t)}\left(\delta_{A} A(t)+\delta_{I} I(t)-(\gamma+\mu) R(t)\right) \\
=& -(\frac{\Lambda}{S}+\gamma\frac{R}{S}+\frac{\beta_{I} IS}{(1+b I)A}+\frac{\beta_{A} S }{(1+b A)}+\alpha\frac{A}{I}+\delta_{A}\frac{A}{R}+\delta_{I}\frac{I}{R})\\
&+m+\frac{\beta_{I} I(t)}{1+b I(t)}+\frac{\beta_{A} A(t)}{1+b A(t)}\\
<&-3\sqrt[3]{\Lambda\beta\alpha}+\beta_{I} I(t)+(\beta_{A}+b\alpha) A(t)\\
&+(4\mu+\alpha+\delta_{A}+\delta_{I}+d+\gamma)+\frac{1}{2}(\sigma_{1}^{2}+\sigma_{2}^{2}+\sigma_{3}^{2}+\sigma_{4}^{2})\\
=&-m(R_{0}^{s}-1)+\beta_{I} I(t)+(\beta_{A}+b\alpha) A(t),
\end{aligned}
\end{equation*}
Similarly, it can be derived that
\begin{equation*}
\begin{aligned}\label{29}
L W_{3}=&-(\lambda)(\alpha A-(\mu+d+\delta_{I}) I),\\
L W_{4}  =&-\frac{1}{S(t)}\left(\Lambda -(\frac{\beta_{I} I(t)}{1+b I(t)}+\frac{\beta_{A} A(t)}{1+b A(t)}) S(t)-\mu S(t)+\gamma R(t)\right)\\
&-\frac{1}{I(t)}\left(\alpha A(t)-\left(d+\delta_{I}+\mu\right) I(t)\right)-\frac{1}{R(t)}(\delta_{A}A(t)+\delta_{I}I(t)\\
&-(\mu+\gamma)R(t))+\frac{1}{2} \sigma_{1}^{2}+\frac{1}{2} \sigma_{2}^{2} +\frac{1}{2}\sigma_{4}^{2}\\
=&  -\frac{\Lambda}{S}+\frac{\beta_{I} I(t)}{1+b I(t)}+\frac{\beta_{A} A(t)}{1+bA(t)}-\gamma\frac{R(t)}{S(t)}-\frac{\alpha A(t)}{I(t)}-\frac{\delta_{A}A(t)}{R(t)}\\
&-\frac{\delta_{I}I(t)}{R(t)}+\left(d+\delta_{I}+3\mu+\gamma\right)+\frac{1}{2}\sigma_{1}^{2}+\frac{1}{2} \sigma_{2}^{2}+\frac{1}{2} \sigma_{4}^{2}.
\end{aligned}
\end{equation*}
Let  $\lambda=\frac{\beta_{A}+b\alpha}{\alpha}$. Thus,
\begin{eqnarray}\label{31}
L W& <& H-\frac{B}{2}(S^{1+\theta}+A^{1+\theta}+I^{1+\theta}+R^{1+\theta})+(-m(R_{0}^{s}-1)+(\beta_{I}+\frac{(\beta_{A}+b\alpha)(d+\mu+\delta_{I})}{\alpha}) I\nonumber\\
&&-\frac{\Lambda}{S}+\frac{\beta_{I} I(t)}{1+b I(t)}+\frac{\beta_{A} A(t)}{1+bA(t)}-\gamma\frac{R(t)}{S(t)}-\frac{\alpha A(t)}{I(t)}-\frac{\delta_{A}A(t)}{R(t)}\nonumber\\
&&-\frac{\delta_{I}I(t)}{R(t)}+\left(d+\delta_{I}+3\mu+\gamma\right)+\frac{1}{2}\sigma_{1}^{2}+\frac{1}{2} \sigma_{2}^{2}+\frac{1}{2} \sigma_{4}^{2}.
\end{eqnarray}

Following Lemma \ref{L2},  a compact set $U_{\epsilon}$ is described as
$$
U_{\epsilon}=\left\{(S(t), A(t), I(t), R(t)) \in \mathbb{R}_{+}^{4} \mid \epsilon \leq S(t) \leq \frac{1}{\epsilon}, \epsilon \leq A(t) \leq \frac{1}{\epsilon}, \epsilon^{2} \leq I(t) \leq \frac{1}{\epsilon^{2}}, \epsilon^{3} \leq R(t) \leq \frac{1}{\epsilon^{3}}\right\},
$$
where $\epsilon>0 $ is a sufficiently small constant and satisfies the following inequality:
\begin{align}
\label{32}
&-\frac{\Lambda}{\epsilon}+F_{1} \leq-1, \\ \label{33}
&-\alpha \frac{1}{\epsilon}+F_{2} \leq-1,\\ \label{34}
&-\delta_{A} \frac{1}{\epsilon^{2}}+F_{3} \leq-1,\\ \label{35}
&-\delta_{I} \frac{1}{\epsilon}+F_{4} \leq-1,\\ \label{36}
&-\frac{B}{4}\left(\frac{1}{\epsilon}\right)^{\theta+1}+F_{5} \leq-1,\\ \label{37}
&-\frac{B}{4}\left(\frac{1}{\epsilon}\right)^{\theta+1}+F_{6} \leq-1,\\ \label{38}
&-\frac{B}{4}\left(\frac{1}{\epsilon}\right)^{2(\theta+1)}+F_{7} \leq-1,\\ \label{39}
&-\frac{B}{4}\left(\frac{1}{\epsilon}\right)^{3(\theta+1)}+F_{8} \leq-1,
\end{align}
where  $F_{i}(i=1,2, \ldots, 8)$  are determined by formula (\ref{31}).  Further,  the domain  $\mathbb{R}_{+}^{4} \backslash U_{\epsilon}$ is divided into eight domains:
\begin{equation}\nonumber
\begin{array}{ll}
U_{1}=\left\{(S(t), A(t), I(t),R(t)) \in \mathbb{R}_{+}^{4}| S(t)<\epsilon\right\},\\
U_{2}=\left\{(S(t), A(t), I(t),R(t)) \in \mathbb{R}_{+}^{4}| I(t)<\epsilon^{2},A(t)>\epsilon\right\}, \\
U_{3}=\left\{(S(t), A(t), I(t),R(t)) \in \mathbb{R}_{+}^{4}| R(t)<\epsilon^{3},A(t)>\epsilon\right\}, \\
U_{4}=\left\{(S(t), A(t), I(t),R(t))\in  \mathbb{R}_{+}^{4}| R(t)<\epsilon^{3},I(t)>\epsilon^{2}\right\} ,\\
U_{5}=\left\{(S(t), A(t), I(t),R(t))\in \mathbb{R}_{+}^{4}| S(t)>\frac{1}{\epsilon}\right\}, U_{6}=\left\{(S(t), A(t), I(t),R(t)) \in \mathbb{R}_{+}^{4}| A(t)>\frac{1}{\epsilon}\right\},\\
U_{7}=\left\{(S(t), A(t), I(t),R(t)) \in \mathbb{R}_{+}^{4}|I(t)>\frac{1}{\epsilon^{2}}\right\}, U_{8}=\left\{(S(t), A(t), I(t),R(t)) \in \mathbb{R}_{+}^{4}| R(t)>\frac{1}{\epsilon^{3}}\right\}.
\end{array}
\end{equation}
Next, we will prove $LW<-1$ on each region $U_i(i=1,2,\ldots,7)$.
\begin{case}
 \textnormal{ When $ (S, A, I,R) \in U_{1}$, according to (\ref{31}) and (\ref{32}), it follows that
\begin{equation}\nonumber
\begin{aligned}
  L W &<-\frac{\Lambda}{S}+H-\frac{B}{2}(S^{1+\theta}+A^{1+\theta}+I^{1+\theta}+R^{1+\theta})+(-m(R_{0}^{s}-1)+(\beta_{I}+\frac{(\beta_{A}+b\alpha)(d+\mu+\delta_{I})}{\alpha}) I\\
&+\frac{\beta_{I} I(t)}{1+b I(t)}+\frac{\beta_{A} A(t)}{1+bA(t)}-\gamma\frac{R(t)}{S(t)}-\frac{\alpha A(t)}{I(t)}-\frac{\delta_{A}A(t)}{R(t)}-\frac{\delta_{I}I(t)}{R(t)}+\left(d+\delta_{I}+3\mu+\gamma\right)\\
&+\frac{1}{2}\sigma_{1}^{2}+\frac{1}{2} \sigma_{2}^{2}+\frac{1}{2} \sigma_{4}^{2}\leq-\frac{\Lambda}{\epsilon}+F_{1} \leq-1,
\end{aligned}
\end{equation}
where
\begin{equation}\nonumber
\begin{aligned}
F_{1}= & \sup _{(S(t), A(t), I(t), R(t)) \in \mathbb{R}_{+}^{4}}\{H-\frac{B}{2}(S^{1+\theta}+A^{1+\theta}+R^{1+\theta})+\frac{\beta_{I}+\beta_{A}}{b}\\
&-m(R_{0}^{s}-1)-(\frac{B}{2}I^{\theta}-\frac{\alpha\beta_{I}+(\beta_{A}+b\alpha)(d+\mu+\delta_{I})}{\alpha})I-\gamma\frac{R(t)}{S(t)}-\frac{\alpha A(t)}{I(t)}\\
&-\frac{\delta_{A}A(t)}{R(t)}-\frac{\delta_{I}I(t)}{R(t)}+\left(d+\delta_{I}+3\mu+\gamma\right)+\frac{1}{2}\sigma_{1}^{2}+\frac{1}{2} \sigma_{2}^{2}+\frac{1}{2} \sigma_{4}^{2}\}.
\end{aligned}
\end{equation}
}
\end{case}
\begin{case}
\textnormal{ When $ (S, A, I,R) \in U_{2}$, based on (\ref{31}) and (\ref{33}), one has
\begin{equation}\nonumber
\begin{aligned}
  L W <&-\frac{\alpha A(t)}{I(t)}+H-\frac{B}{2}(S^{1+\theta}+A^{1+\theta}+I^{1+\theta}+R^{1+\theta})+(\beta_{I}+\frac{(\beta_{A}+b\alpha)(d+\mu+\delta_{I})}{\alpha}) I\\
&-\frac{\Lambda}{S(t)}+\frac{\beta_{I} I(t)}{1+b I(t)}+\frac{\beta_{A} A(t)}{1+bA(t)}-\gamma\frac{R(t)}{S(t)}-\frac{\delta_{A}A(t)}{R(t)}-\frac{\delta_{I}I(t)}{R(t)}+\left(d+\delta_{I}+3\mu+\gamma\right)\\
&+\frac{1}{2}\sigma_{1}^{2}+\frac{1}{2} \sigma_{2}^{2}+\frac{1}{2} \sigma_{4}^{2}\leq-\frac{\alpha}{\epsilon}+F_{2} \leq-1,
\end{aligned}
\end{equation}
where
\begin{equation}\nonumber
\begin{aligned}
F_{2}= & \sup _{(S(t), A(t), I(t), R(t)) \in \mathbb{R}_{+}^{4}}\{H-\frac{B}{2}(S^{1+\theta}+A^{1+\theta}+R^{1+\theta})-\frac{\Lambda}{S(t)}+\frac{\beta_{I}+\beta_{A}}{b}\\
&-m(R_{0}^{s}-1)-(\frac{B}{2}I^{\theta}-\frac{\alpha\beta_{I}+(\beta_{A}+b\alpha)(d+\mu+\delta_{I})}{\alpha})I-\gamma\frac{R(t)}{S(t)}-\frac{\delta_{A}A(t)}{R(t)}-\frac{\delta_{I}I(t)}{R(t)}\\
&+\left(d+\delta_{I}+3\mu+\gamma\right)+\frac{1}{2}\sigma_{1}^{2}+\frac{1}{2} \sigma_{2}^{2}+\frac{1}{2} \sigma_{4}^{2}\}.
\end{aligned}
\end{equation}
}
\end{case}
\begin{case}
\textnormal{ When $ (S, A, I,R) \in U_{3}$, according to (\ref{31}) and (\ref{34}), it follows that
\begin{equation}\nonumber
\begin{aligned}
  L W <&-\frac{\delta_{A}A(t)}{R(t)}+H-\frac{B}{2}(S^{1+\theta}+A^{1+\theta}+I^{1+\theta}+R^{1+\theta})+(\beta_{I}+\frac{(\beta_{A}+b\alpha)(d+\mu+\delta_{I})}{\alpha}) I\\
&-\frac{\Lambda}{S(t)}+\frac{\beta_{I} I(t)}{1+b I(t)}+\frac{\beta_{A} A(t)}{1+bA(t)}-\frac{\alpha A(t)}{I(t)}-\gamma\frac{R(t)}{S(t)}-\frac{\delta_{I}I(t)}{R(t)}+\left(d+\delta_{I}+3\mu+\gamma\right)\\
&+\frac{1}{2}\sigma_{1}^{2}+\frac{1}{2} \sigma_{2}^{2}+\frac{1}{2} \sigma_{4}^{2}\leq-\frac{\delta_{A}}{\epsilon^{2}}+F_{3} \leq-1,
\end{aligned}
\end{equation}
where
\begin{equation}\nonumber
\begin{aligned}
F_{3}= & \sup _{(S(t), A(t), I(t), R(t)) \in \mathbb{R}_{+}^{4}}\{H-\frac{B}{2}(S^{1+\theta}+A^{1+\theta}+R^{1+\theta})-\frac{\Lambda}{S(t)}+\frac{\beta_{I}+\beta_{A}}{b}\\
&-m(R_{0}^{s}-1)-(\frac{B}{2}I^{\theta}-\frac{\alpha\beta_{I}+(\beta_{A}+b\alpha)(d+\mu+\delta_{I})}{\alpha})I-\gamma\frac{R(t)}{S(t)}-\frac{\alpha A(t)}{I(t)}-\frac{\delta_{I}I(t)}{R(t)}\\
&+\left(d+\delta_{I}+3\mu+\gamma\right)+\frac{1}{2}\sigma_{1}^{2}+\frac{1}{2} \sigma_{2}^{2}+\frac{1}{2} \sigma_{4}^{2}\}.
\end{aligned}
\end{equation}
}
\end{case}
\begin{case}
\textnormal{ When $ (S, A, I,R) \in U_{4}$, combining (\ref{31}) and (\ref{35}) yields
\begin{equation}\nonumber
\begin{aligned}
  L W <&-\frac{\delta_{I}I(t)}{R(t)}+H-\frac{B}{2}(S^{1+\theta}+A^{1+\theta}+I^{1+\theta}+R^{1+\theta})+(-m(R_{0}^{s}-1)+(\beta_{I}+\frac{(\beta_{A}+b\alpha)(d+\mu+\delta_{I})}{\alpha})I\\
&-\frac{\Lambda}{S(t)}+\frac{\beta_{I} I(t)}{1+b I(t)}+\frac{\beta_{A} A(t)}{1+bA(t)}-\frac{\alpha A(t)}{I(t)}-\gamma\frac{R(t)}{S(t)}-\frac{\delta_{A}A(t)}{R(t)}+\left(d+\delta_{I}+3\mu+\gamma\right)\\
&+\frac{1}{2}\sigma_{1}^{2}+\frac{1}{2} \sigma_{2}^{2}+\frac{1}{2} \sigma_{4}^{2}\leq-\frac{\delta_{I}}{\epsilon}+F_{4} \leq-1,
\end{aligned}
\end{equation}
where
\begin{equation}\nonumber
\begin{aligned}
F_{4}= & \sup _{(S(t), A(t), I(t), R(t)) \in \mathbb{R}_{+}^{4}}\{H-\frac{B}{2}(S^{1+\theta}+A^{1+\theta}+R^{1+\theta})-\frac{\Lambda}{S(t)}+\frac{\beta_{A}+\beta_{I}}{b}\\
&-m(R_{0}^{s}-1)-(\frac{B}{2}I^{\theta}-\frac{\alpha\beta_{I}+(\beta_{A}+b\alpha)(d+\mu+\delta_{I})}{\alpha})I-\gamma\frac{R(t)}{S(t)}-\frac{\alpha A(t)}{I(t)}-\frac{\delta_{A}A(t)}{R(t)}\\
&+\left(d+\delta_{I}+3\mu+\gamma\right)+\frac{1}{2}\sigma_{1}^{2}+\frac{1}{2} \sigma_{2}^{2}+\frac{1}{2} \sigma_{4}^{2}\}.
\end{aligned}
\end{equation}}
\end{case}
\begin{case}
\textnormal{ When $ (S, A, I,R) \in U_{5}$, (\ref{31}) and (\ref{36}) imply that
\begin{equation}\nonumber
\begin{aligned}
  L W <&-\frac{B}{2}S^{1+\theta}+H-\frac{B}{2}(A^{1+\theta}+I^{1+\theta}+R^{1+\theta})+(-m(R_{0}^{s}-1)+(\beta_{I}+\frac{(\beta_{A}+b\alpha)(d+\mu+\delta_{I})}{\alpha})I\\
&-\frac{\Lambda}{S(t)}+\frac{\beta_{I} I(t)}{1+b I(t)}+\frac{\beta_{A} A(t)}{1+bA(t)}-\frac{\alpha A(t)}{I(t)}-\gamma\frac{R(t)}{S(t)}-\frac{\delta_{A}A(t)}{R(t)}-\frac{\delta_{I}I(t)}{R(t)}\\
&+\left(d+\delta_{I}+3\mu+\gamma\right)+\frac{1}{2}\sigma_{1}^{2}+\frac{1}{2} \sigma_{2}^{2}+\frac{1}{2} \sigma_{4}^{2}\leq-\frac{B}{4}(\frac{1}{\epsilon})^{\theta+1}+F_{5} \leq-1,
\end{aligned}
\end{equation}
where
\begin{equation}\nonumber
\begin{aligned}
F_{5}= & \sup _{(S(t), A(t), I(t), R(t)) \in \mathbb{R}_{+}^{4}}\{-\frac{B}{4}S^{1+\theta}+H-\frac{B}{2}(A^{1+\theta}+R^{1+\theta}) -\frac{\Lambda}{S(t)}+\frac{\beta_{A}+\beta_{I}}{b}\\
&-m(R_{0}^{s}-1)-(\frac{B}{2}I^{\theta}-\frac{\alpha\beta_{I}+(\beta_{A}+b\alpha)(d+\mu+\delta_{I})}{\alpha})I-\gamma\frac{R(t)}{S(t)}-\frac{\alpha A(t)}{I(t)}-\frac{\delta_{A}A(t)}{R(t)}-\frac{\delta_{I}I(t)}{R(t)}\\
&+\left(d+\delta_{I}+3\mu+\gamma\right)+\frac{1}{2}\sigma_{1}^{2}+\frac{1}{2} \sigma_{2}^{2}+\frac{1}{2} \sigma_{4}^{2}\}.
\end{aligned}
\end{equation}
}
\end{case}
\begin{case}
\textnormal{ When $ (S, A, I,R) \in U_{6}$, based on(\ref{31}) and (\ref{37}), it follows that
\begin{equation}\nonumber
\begin{aligned}
  L W<&-\frac{B}{2}A^{1+\theta}+H-\frac{B}{2}(S^{1+\theta}+I^{1+\theta}+R^{1+\theta})+(-m(R_{0}^{s}-1)+(\beta_{I}+\frac{(\beta_{A}+b\alpha)(d+\mu+\delta_{I})}{\alpha}) I\\
&-\frac{\Lambda}{S(t)}+\frac{\beta_{I} I(t)}{1+b I(t)}+\frac{\beta_{A} A(t)}{1+bA(t)}-\frac{\alpha A(t)}{I(t)}-\gamma\frac{R(t)}{S(t)}-\frac{\delta_{A}A(t)}{R(t)}-\frac{\delta_{I}I(t)}{R(t)}\\
&+\left(d+\delta_{I}+3\mu+\gamma\right)+\frac{1}{2}\sigma_{1}^{2}+\frac{1}{2} \sigma_{2}^{2}+\frac{1}{2} \sigma_{4}^{2}\leq-\frac{B}{4}(\frac{1}{\epsilon})^{\theta+1}+F_{6} \leq-1,
\end{aligned}
\end{equation}
where
\begin{equation}\nonumber
\begin{aligned}
F_{6}= & \sup _{(S(t), A(t), I(t), R(t)) \in \mathbb{R}_{+}^{4}}\{-\frac{B}{4}A^{1+\theta}+H-\frac{B}{2}(S^{1+\theta}+R^{1+\theta}) -\frac{\Lambda}{S(t)}+\frac{\beta_{A}+\beta_{I}}{b}\\
&-m(R_{0}^{s}-1)-(\frac{B}{2}I^{\theta}-\frac{\alpha\beta_{I}+(\beta_{A}+b\alpha)(d+\mu+\delta_{I})}{\alpha})I-\gamma\frac{R(t)}{S(t)}-\frac{\alpha A(t)}{I(t)}-\frac{\delta_{A}A(t)}{R(t)}-\frac{\delta_{I}I(t)}{R(t)}\\
&+\left(d+\delta_{I}+3\mu+\gamma\right)+\frac{1}{2}\sigma_{1}^{2}+\frac{1}{2} \sigma_{2}^{2}+\frac{1}{2} \sigma_{4}^{2}\}.
\end{aligned}
\end{equation}
}
\end{case}
\begin{case}
\textnormal{ When $ (S, A, I,R) \in U_{7}$, (\ref{31}) and (\ref{38}) imply that
\begin{equation}\nonumber
\begin{aligned}
  L W<&-\frac{B}{2}I^{1+\theta}+H-\frac{B}{2}(S^{1+\theta}+A^{1+\theta}+R^{1+\theta})+(-m(R_{0}^{s}-1)+(\beta_{I}+\frac{(\beta_{A}+b\alpha)(d+\mu+\delta_{I})}{\alpha}) I\\
&-\frac{\Lambda}{S(t)}+\frac{\beta_{I} I(t)}{1+b I(t)}+\frac{\beta_{A} A(t)}{1+bA(t)}-\frac{\alpha A(t)}{I(t)}-\gamma\frac{R(t)}{S(t)}-\frac{\delta_{A}A(t)}{R(t)}-\frac{\delta_{I}I(t)}{R(t)}\\
&+\left(d+\delta_{I}+3\mu+\gamma\right)+\frac{1}{2}\sigma_{1}^{2}+\frac{1}{2} \sigma_{2}^{2}+\frac{1}{2} \sigma_{4}^{2}\leq-\frac{B}{4}(\frac{1}{\epsilon})^{2(\theta+1)}+F_{7} \leq-1,
\end{aligned}
\end{equation}
where
\begin{equation}\nonumber
\begin{aligned}
F_{7}= & \sup _{(S(t), A(t), I(t), R(t)) \in \mathbb{R}_{+}^{4}}\{-(\frac{B}{4}I^{\theta}-\frac{\alpha\beta_{I}+(\beta_{A}+b\alpha)(d+\mu+\delta_{I})}{\alpha})I-\frac{B}{2}(S^{1+\theta}+A^{1+\theta}+R^{1+\theta})\\
&+H-\frac{\Lambda}{S(t)}+\frac{\beta_{A}+\beta_{I}}{b}-m(R_{0}^{s}-1)-\gamma\frac{R(t)}{S(t)}-\frac{\alpha A(t)}{I(t)}-\frac{\delta_{A}A(t)}{R(t)}-\frac{\delta_{I}I(t)}{R(t)}\\
&+\left(d+\delta_{I}+3\mu+\gamma\right)+\frac{1}{2}\sigma_{1}^{2}+\frac{1}{2} \sigma_{2}^{2}+\frac{1}{2} \sigma_{4}^{2}\}.
\end{aligned}
\end{equation}
}
\end{case}
\begin{case}
\textnormal{ When $ (S, A, I,R) \in U_{8}$,   from (\ref{31}) and (\ref{39}), one has
\begin{equation}\nonumber
\begin{aligned}
  L W < &-\frac{B}{2}R^{1+\theta}+H-\frac{B}{2}(S^{1+\theta}+A^{1+\theta}+I^{1+\theta})+(-m(R_{0}^{s}-1)+(\beta_{I}+\frac{(\beta_{A}+b\alpha)(d+\mu+\delta_{I})}{\alpha}) I\\
&-\frac{\Lambda}{S(t)}+\frac{\beta_{I} I(t)}{1+b I(t)}+\frac{\beta_{A} A(t)}{1+bA(t)}-\frac{\alpha A(t)}{I(t)}-\gamma\frac{R(t)}{S(t)}-\frac{\delta_{A}A(t)}{R(t)}-\frac{\delta_{I}I(t)}{R(t)}\\&+\left(d+\delta_{I}+3\mu+\gamma\right)+\frac{1}{2}\sigma_{1}^{2}+\frac{1}{2} \sigma_{2}^{2}+\frac{1}{2} \sigma_{4}^{2}\leq-\frac{B}{4}(\frac{1}{\epsilon})^{3(\theta+1)}+F_{8} \leq-1,
\end{aligned}
\end{equation}
where
\begin{equation}\nonumber
\begin{aligned}
F_{8}= & \sup _{(S(t), A(t), I(t), R(t)) \in \mathbb{R}_{+}^{4}}\{-\frac{B}{4}R^{1+\theta}+H-\frac{B}{2}(S^{1+\theta}+A^{1+\theta}) -\frac{\Lambda}{S(t)}+\frac{\beta_{A}+\beta_{I}}{b}\\
&-m(R_{0}^{s}-1)-(\frac{B}{2}I^{\theta}-\frac{\alpha\beta_{I}+(\beta_{A}+b\alpha)(d+\mu+\delta_{I})}{\alpha})I-\gamma\frac{R(t)}{S(t)}-\frac{\alpha A(t)}{I(t)}-\frac{\delta_{A}A(t)}{R(t)}-\frac{\delta_{I}I(t)}{R(t)}\\
&+\left(d+\delta_{I}+3\mu+\gamma\right)+\frac{1}{2}\sigma_{1}^{2}+\frac{1}{2} \sigma_{2}^{2}+\frac{1}{2} \sigma_{4}^{2}\}.
\end{aligned}
\end{equation}
}
\end{case}
In conclusion, $LW < -1$ always holds in the eight regions of the set $U$.
Moreover,   the diffusion matrix of system (\ref{eq2}) is expressed by
\begin{equation}\nonumber
C_{4 \times 4}=\left[\begin{array}{cccc}
\sigma_{1}^{2} S^{2} & 0 & 0 & 0 \\
0 & \sigma_{2}^{2} A^{2} & 0 & 0 \\
0 & 0 & \sigma_{3}^{2} I^{2}& 0 \\
0 & 0 &0 & \sigma_{4}^{2} R^{2}
\end{array}\right].
\end{equation}
Let  $M_{0}:=\min _{\left(S, A, I, R\right) \in \mathscr{H}_{k} \subset \mathbb{R}_{+}^{4}}\left\{\sigma_{1}^{2} S^{2}, \sigma_{2}^{2} A^{2}, \sigma_{3}^{2} I^{2}, \sigma_{4}^{2} R^{2}\right\}$, for any $ k>1$,
$\mathscr{H}_{k}:= [1/k, k ] \times [1/k, k] \times [1/k^{2}, k^{2}] \times[1/k^{3}, k^{3}]$. Then
\begin{equation*}
\begin{aligned}
\sum_{i, j=1}^{4} a_{i j}\left(S, A, I, R\right) \eta_{i} \eta_{j} & = (
\sigma_{1} S\eta_{1} \ \ \sigma_{2} A \eta_{2}\ \  \sigma_{3} I \eta_{3} \ \ \sigma_{4} R\eta_{4})
(\sigma_{1} S \eta_{1} \ \
\sigma_{2} A \eta_{2} \ \
\sigma_{3} I \eta_{3} \ \
\sigma_{4} R \eta_{4})' \\
& =\left(\sigma_{1} S\right)^{2} \eta_{1}^{2}+\left(\sigma_{2} A\right)^{2} \eta_{2}{ }^{2}+\left(\sigma_{3} I\right)^{2} \eta_{3}{ }^{2}+\left(\sigma_{4} R\right)^{2} \eta_{4}{ }^{2} \geq M_{0}\left\|\eta^{2}\right\|
\end{aligned}
\end{equation*}
for all  $\left(S,A, I,R\right) \in \mathscr{H}_{k}$  and $ \eta=\left(\eta_{1}, \eta_{2}, \eta_{3}, \eta_{4}\right) \in \mathbb{R}^{4}$.
In summary, system (\ref{eq2}) has a unique ergodic stationary distribution.
\end{proof}

\section{Stochastic optimal control}
After investigating the dynamical properties, we will study stochastic optimal control of the SAIRS epidemic model. The study is carried out with vaccination of susceptible individuals and isolation of asymptomatic infected ones as control variables, denoted by $u_{1}(t)$ and   $u_{2}(t)$, in the SAIRS model (\ref{eq2}).  The effect of the control variable $u_{1}(t)$  is a reduction in the number of susceptible people, and that of the variable $u_{2}(t)$ is to decrease the transmission rate. Thus, we establish the following stochastic controlled system
\begin{equation}\label{eq3}
\left\{\begin{array}{l}
		d S(t)=\left[\Lambda-\left(\frac{\beta_{I} I(t)}{1+b I(t)}+\frac{\beta_{A} A(t)}{1+b A(t)}\right) (1-u_{2}(t))S(t)-(\mu+u_{1}(t)) S(t)+\gamma R(t)\right]dt+\sigma_{1} S(t) d B_{1}(t),\\
		d A(t)=\left[\left(\frac{\beta_{I} I(t)}{1+b I(t)}+\frac{\beta_{A} A(t)}{1+b A(t)}\right) S(t)(1-u_{2}(t)) S(t)-\left(\alpha+\delta_{A}+\mu\right) A(t)\right]dt+\sigma_{2}A(t) d B_{2}(t), \\
		d I(t)=[\alpha A(t)-\left(\delta_{I}+\mu+d\right) I(t) ]dt+\sigma_{3} I(t) dB_{3}(t),\\
		d R(t)=[\delta_{A} A(t)+\delta_{I} I(t)+ u_{1}(t)S(t)-(\gamma+\mu) R(t)]dt+\sigma_{4}R(t) d B_{4}(t),
	\end{array}
\right.
\end{equation}
where initial values
$S(0)>0$, $A(0) \geqslant 0$, $I(0) \geqslant 0$, $R(0)>0.$
For convenience, denote   $u(t)$=$(u_{1}(t),$  $u_{2}(t))^{\prime}$, $x(t)$ $=(x_{1}(t)$, $x_{2}(t)$, $x_{3}(t)$, $x_{4}(t))^{\prime}=(S(t),A(t), I(t), R(t))^{\prime}$, and
\begin{equation}\label{eq8}\nonumber
\begin{aligned}
f_{1}(x(t),u(t))=&\left[\Lambda-\left(\frac{\beta_{I} I(t)}{1+b I(t)}+\frac{\beta_{A} A(t)}{1+b A(t)}\right) (1-u_{2}(t)) S(t)-(\mu+u_{1}(t)) S(t) +\gamma R(t)\right]dt, \\	
f_{2}(x(t),u(t))=&\left[\left(\frac{\beta_{I} I(t)}{1+b I(t)}+\frac{\beta_{A} A(t)}{1+b A(t)}\right)(1-u_{2}(t)) S(t)-\left(\alpha+\delta_{A}+\mu\right) A(t)\right]dt\\
f_{3}(x(t),u(t))=&[\alpha A(t)-\left(\delta_{I}+\mu+d\right) I(t)]dt,\\
f_{4}(x(t),u(t))=&[\delta_{A} A(t)+\delta_{I} I(t)+ u_{1}(t)S(t)-(\gamma+\mu) R(t)]dt,
	\end{aligned}
\end{equation}	
and  $g_{1}(x(t))=\sigma_{1} S(t), g_{2}(x(t))=\sigma_{2} A(t), g_{3}(x(t))=\sigma_{3} I(t), g_{4}(x(t))=\sigma_{4} R(t).$ Then,  the model (\ref{eq3}) is equivalent to
$$d x(t)=f(x(t), u(t)) d t+g(x(t)) d B(t)$$ with initial value
$x(0)=(x_{1}(0), x_{2}(0), x_{3}(0), x_{4}(0))^{\prime}=x(0),$
where  $f=(f_1,f_2,f_3,f_4)^{\prime}$, $g=(g_1,g_2,g_3,g_4)^{\prime}$, and $B(t)=(B_1(t),B_2(t),B_3(t),B_4(t))^{\prime}$. We are interested in the objective function as
\begin{equation*}\label{eq10}
	\begin{aligned}
		J(u)=&\frac{E}{2}\left\{\int_{0}^{T}\left(P_{1} S+P_{2} A+P_{3} I+\frac{Q_{1}}{2} u_{1}^{2}(t)+\frac{Q_{2}}{2} u_{2}^{2}(t)\right) d t+\frac{k_{1}}{2} S^{2}+\frac{k_{2}}{2} A^{2}+\frac{k_{3}}{2} I^{2}+\frac{k_{4}}{2} R^{2}\right\},
	\end{aligned}
\end{equation*}
where $ P_{i}(i=1,2,3), Q_{i}(i=1,2), k_{i}(i=1,2,3,4)$  are  non-negative constants. We aim to obtain an optimal strategy $u^{*}(t)=(u_{1}^{*}(t), u_{2}^{*}(t))$ such that
$J(u) \geqslant J\left(u^{*}\right)$ for any $ u \in U$, where
$U$ is a controlling admissible set and
\begin{equation}\label{12}\nonumber
  \begin{aligned}
U&:=\{\left(u_{1}, u_{2}\right)| 0 \leqslant u_{i} \leqslant 1,  u_{i} \text {  is Lebesgue measurable on }[0, T] ,  i =1,2\}.
\end{aligned}
\end{equation}
\begin{theorem}
 The controlled  system (\ref{eq3}) has an optimal control $u^{*}(t)=(u_{1}^{*}(t), u_{2}^{*}(t))$  satisfying  $J(u) \geqslant J\left(u^{*}\right)$ for any $ u \in U$.
\end{theorem}
\begin{proof}
Through the Pontryagin's Minimum Principle \textcolor{blue}{(Yong et al. 1999)},  the Hamiltonian function $H(x, u, m, n)$  is formed by
\begin{equation}\label{eq13}\nonumber
H(x, u, m, n) = \langle f(x, u), m\rangle+l(x, u)+\langle g(x), n\rangle,
\end{equation}
where $ \langle. , .\rangle $ means the Euclidean inner product and $ m=(m_{1}, m_{2}, m_{3}, m_{4})^{\prime} $ and  $n=(n_{1}, n_{2}, n_{3}, n_{4})^{\prime}  $ are two adjoint vectors. Based on the minimum principle, we obtain
\begin{equation*}\label{eq14}\nonumber
	d x^{*}(t)=\frac{\partial H\left(x^{*}(t), u^{*}(t), m, n\right)}{\partial m} d t+g\left(x^{*}(t)\right) d B(t),\quad
	d m^{*}(t)=-\frac{\partial H\left(x^{*}(t), u^{*}(t), m, n\right)}{\partial x} d t+n(t) d B(t),
\end{equation*}
and
	\begin{equation}\label{eq16}
	H_{m}\left(x^{*}(t), u^{*}(t), m, n\right)=\min _{u \in \mathrm{U}} H_{m}\left(x^{*}(t), u^{*}(t), m, n\right),
\end{equation}
where $ x^{*}(t)=(S^*(t),A^*(t), I^*(t), R^*(t))' $ is the path of optimality of  $x(t)$. The initial  states satisfy $x_{0}=x^{*}(0)$ and termination satisfy
$-\frac{\partial h(x^{*}(t_{f}))}{\partial x}=m(t_{f})$, where $h(S, A, I, R)=\frac{k_{1}}{2} S^{2}+\frac{k_{2}}{2} A^{2}+\frac{k_{3}}{2} I^{2}+\frac{k_{4}}{2} R^{2}$ are terminal conditions.  Since (\ref{eq14}) shows that the optimal variable $ x^{*}(t)$ is a vector consisting of $n(t)$, $m(t)$ and $x^{*}(t)$, this implies that
$u^{*}(t)=\phi\left(m, n, x^{*}\right) \text {.}$ $u^{*}(t)$  can be calculated by equation (\ref{eq16}). Thus, the  Hamiltonian function is
\begin{equation*}\label{eq20}
    \begin{aligned}
		H=&P_{1} S+P_{2} A+P_{3} I+u_{1}^{2} \frac{Q_{1}}{2}+u_{2}^{2}\frac{Q_{2}}{2}+S^{2} \frac{k_{1}}{2}+A^{2} \frac{k_{3}}{2}+I^{2} \frac{k_{2}}{2}+R^{2}\frac{k_{4}}{2}\\		&+m_{1}f_{1}(x(t),u(t))+m_{2}f_{2}(x(t),u(t))+m_{3}f_{3}(x(t),u(t))+m_{4}f_{4}(x(t),u(t))\\
		&+\sigma_{1} S n_{1}(t)+\sigma_{2} A n_{2}(t)+\sigma_{3} I n_{3}(t)+\sigma_{4} R n_{4}(t).
	\end{aligned}
\end{equation*}
Using the Pontryagin's Minimum  Principle, we have
\begin{equation*}\label{eq22}
	\begin{aligned}
	   d m_{1}(t)=&\{-P_{1}-\xi_{1}+(m_{1}(t)-m_{2}(t))[(\frac{\beta_{I} I^{*}(t)}{1+b I^{*}(t)}+\frac{\beta_{A} A^{*}(t)}{1+b A^{*}(t)}) (1-u^{*}_{2}(t))]\\
	   &+m_{1}(t)(\mu+u^{*}_{1}(t))-m_{4}(t)u^{*}_{1}(t)-\sigma_{1}n_{1}(t)\}dt+n_{1}(t)dB_{1}(t),\\
	d m_{2}(t)=&\{-P_{2}-\xi_{2}+(m_{1}(t)-m_{2}(t))[\frac{\beta_{A}S^{*}(t)}{(1+bA(t))^{2}}(1-u^{*}_{2}(t))]+m_{2}(t)(\delta_{A}+\mu+\alpha)\\
	&+m_{4}(t)\delta_{A}-\sigma_{2}n_{2}\}dt+n_{2}(t)d B_{2}(t),\\
	d m_{3}(t)=&\{-P_{3}-\xi_{3}+(m_{1}(t)-m_{2})(t)[\frac{\beta_{I}S^{*}(t)}{(1+bI(t))^{2}}(1-u^{*}_{2}(t))]+m_{3}(t)(\delta_{I}+\mu+d)\\
	&+m_{4}(t)\delta_{I}-\sigma_{3}n_{3}(t)\}dt+n_{3}(t)d B_{3}(t),\\
	d m_{4}(t)=&\{-m_{1}(t)\gamma+m_{4}(t)(\gamma+\mu) -\sigma_{4}n_{4}(t)\}dt+n_{4}(t)d B_{4}(t),
	\end{aligned}
\end{equation*}
 where $S^{*}(0)=\widehat{S}, A^{*}(0)=\widehat{A}, B^{*}(0)=\widehat{I}, R^{*}(0)=\widehat{R}$, and $m_{1}\left(T\right)=-k_{1} S$,  $m_{2}\left(T\right)=-A k_{2}$,  $m_{3}\left(T\right)=-I k_{3}$,  $m_{4}\left(T\right)=-R k_{4}.$
The partial derivatives of $ u_{1}$ and $u_{2} $ are calculated by the Hamiltonian function; we obtain the optimal control  $u_{1}^{*}  $ and $  u_{2}^{*}$  as follows
\begin{equation*}\label{eq26}
\begin{array}{l}
u_{1}^{*}=\max \left\{\min \left\{\frac{1}{Q_{1}}\left(m_{1}-m_{4}\right) S^{*}, 1\right\}, 0\right\},
u_{2}^{*}=\max \left\{\min \left\{1, \frac{\left(m_{2}-m_{4}\right) A^{*}+\left(m_{3}-m_{4}\right) I^{*}}{Q_{2}}\right\}, 0\right\}.
\end{array}
\end{equation*}
\end{proof}

\section{ Numerical simulations}
To illustrate the theoretical results,  we perform numerical simulations of the stochastic epidemic system (\ref{eq2}). These simulations can be obtained by  using the stochastic Runge-Kutta method  \textcolor{blue}{(Tocino and Ardanuy 2002)}  as follows
\begin{equation}\nonumber
\begin{aligned}
S_{i+1}=&S_{i}+\left[\Lambda-\left(\frac{\beta_{A} A_{i}}{1+bA_{i}}+\frac{\beta_{I} I_{i}}{1+bI_{i}}\right) S_{i}-\mu S_{i}+\gamma R_{i}\right] \Delta t+\sigma_{1} S_{i} \sqrt{\Delta t} \zeta_{1, i}^{2} +\frac{\sigma_{1}^{2}}{2} S_{i}\left(\zeta_{1, i}^{2}-1\right) \Delta t, \\
A_{i+1}=&A_{i}+\left[\left(\frac{\beta_{A} A_{i}}{1+bA_{i}}+\frac{\beta_{I} I_{i}}{1+bI_{i}}\right) S_{i}-\left(\alpha+\delta_{A}+\mu\right) A_{i}\right]\Delta t+\sigma_{2} A_{i} \sqrt{\Delta t} \zeta_{2, i} +\frac{\sigma_{2}^{2}}{2} A_{i}\left(\zeta_{2, i}^{2}-1\right) \Delta t, \\
I_{i+1}=&I_{i}+\left[\alpha A(t)-\left(\delta_{I}+\mu+d\right) I(t)\right] \Delta t+\sigma_{3} I_{i} \sqrt{\Delta t} \zeta_{3, i} +\frac{\sigma_{3}^{2}}{2} I_{i}\left(\zeta_{3, i}^{2}-1\right) \Delta t, \\
R_{i+1}=&R_{i}+\left[\delta_{A} A(t)+\delta_{I} I(t)-(\gamma+\mu) R(t)\right] \Delta t+\sigma_{4} R_{i} \sqrt{\Delta t} \zeta_{4, i} +\frac{\sigma_{4}^{2}}{2} R_{i}\left(\zeta_{4, i}^{2}-1\right) \Delta t,
\end{aligned}
\end{equation}
where $\zeta_{k, i}(k=1,2,3,4)$ are four independent Gaussian random variables with $N(0,1)$ and the time increment $\Delta t>0$.

\begin{exam} \label{exam1}
Let the initial values $S(0)=1500$, $A(0)=5$, $I(0)=6$, $R(0)=25$, $\Delta t=0.002$, and the parameter values  $\Lambda=30$, $b=0.2$, $\beta_{A}=0.01$,  $\beta_{I}=0.01$, $\gamma=0.5$, $\mu=2\times10^{-5}$, $\delta_{A}=0.2$, $\alpha=0.5$, $\delta_{I}=0.2$,
 $d=0.0027$, $\sigma_{1}=0.05$,  $\sigma_{2}=0.05$, $\sigma_{3}=0.05$, $\sigma_{4}=0.05$ .  Through a simple calculation, $R_{0}^{s}=1.132>1$, satisfying the conditions of Theorem \ref{th2.2}.  Thus, the disease will persist. FIGURE \ref{B0} shows the effects of $b$ from 0.1 and 0.3 on the epidemic system.  It reveals that the spread of the disease decreases as the saturation incidence coefficient increases.
\begin{figure}[!htb]	
\centering	
\begin{minipage}[t]{1\textwidth}
	\centering
	\subfigure[]{
		\includegraphics[width=7.3cm,height=3.9cm]{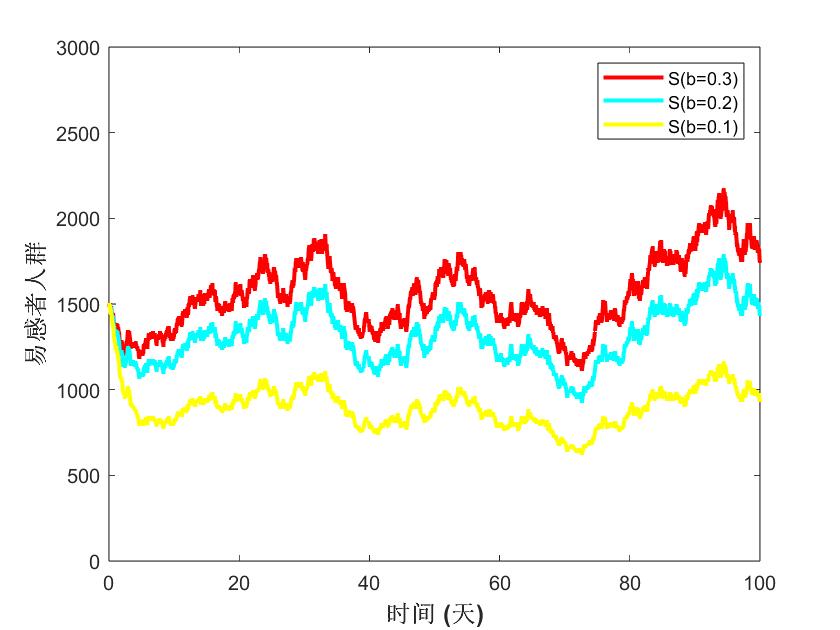} }
	\subfigure[]{
		\includegraphics[width=7.3cm,height=3.9cm]{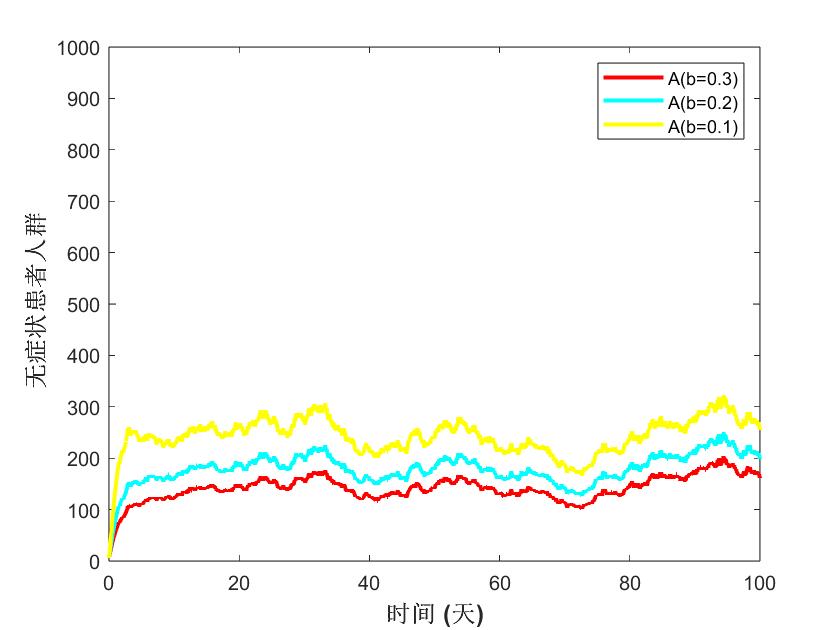} }
	\vspace{-4mm} 
\end{minipage}
\begin{minipage}[t]{1\textwidth}
	\centering
	\subfigure[]{
		\includegraphics[width=7.3cm,height=3.9cm]{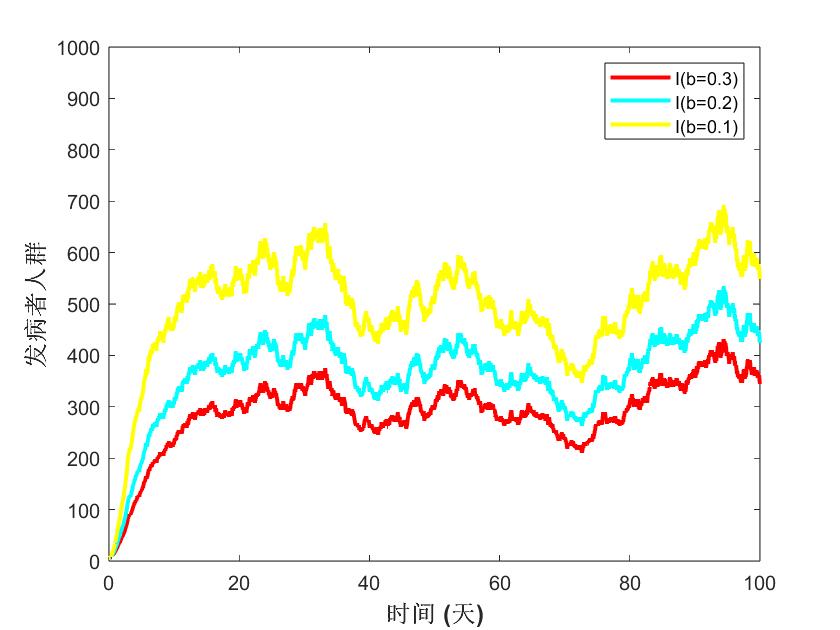} }
	\subfigure[]{
		\includegraphics[width=7.3cm,height=3.9cm]{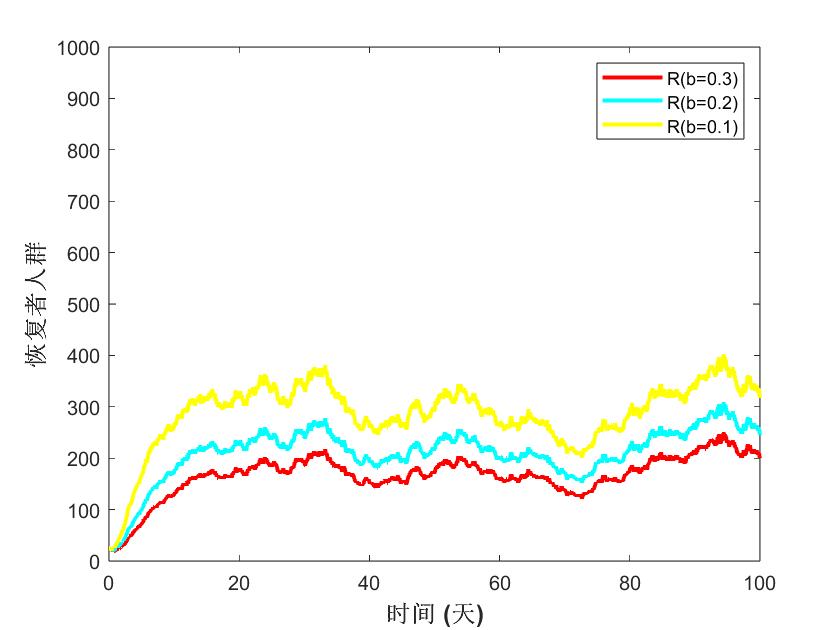} }
\end{minipage}
\vspace{-3mm} 
\caption{ The effect of saturation incidence coefficient $b$ and the persistence of the epidemic system (\ref{eq2}). }
\label{B0}
\end{figure}
\end{exam}

\begin{exam} 
Take $b=0.2$ and  $d=0,0.1,0.2$. Other parameter settings are the same as Example \ref{exam1}.
FIGURE \ref{B1} reflects the impact of fatality rate $d$ on the epidemic system. The number of infected people decreases as the disease-death rate $d$ increases.  The system's various solutions change over time due to differences in fatality rates. From FIGURE \ref{B1},  the greater the diseased-death rate $d$, the smaller the number of people in each compartment.

\begin{figure}[!htb]	
\centering	
\begin{minipage}[t]{1\textwidth}
	\centering
	\subfigure[]{
		\includegraphics[width=7.3cm,height=3.9cm]{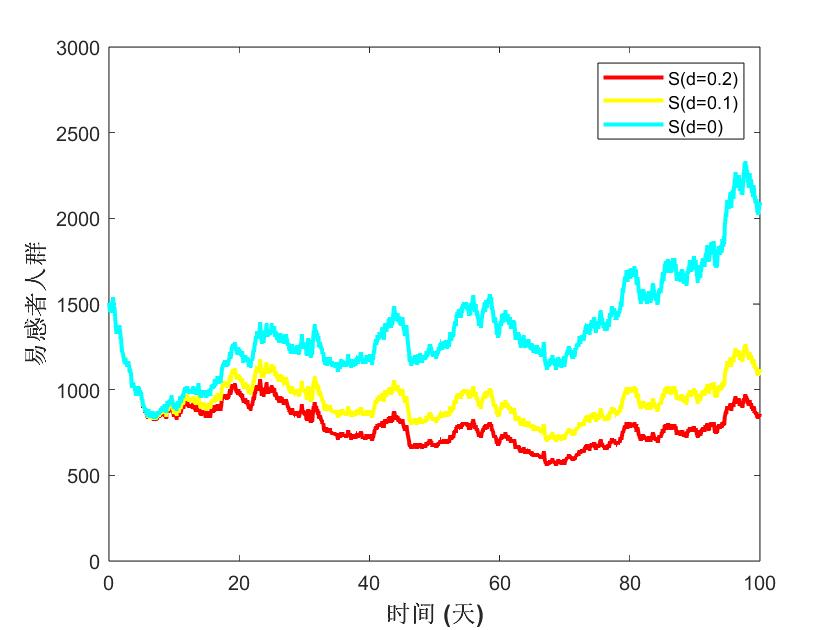} }
	\subfigure[]{
		\includegraphics[width=7.3cm,height=3.9cm]{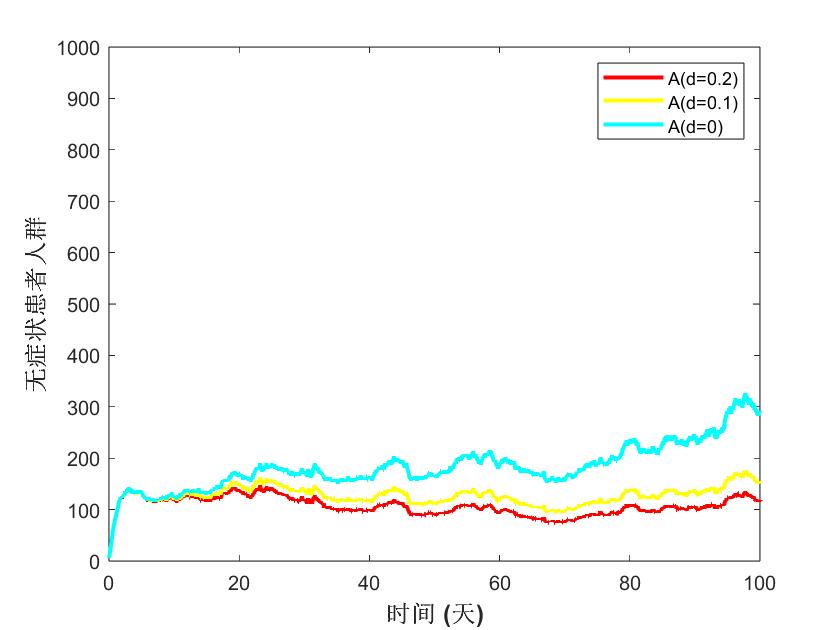} }
	\vspace{-4mm} 
\end{minipage}
\begin{minipage}[t]{1\textwidth}
	\centering
	\subfigure[]{
		\includegraphics[width=7.3cm,height=3.9cm]{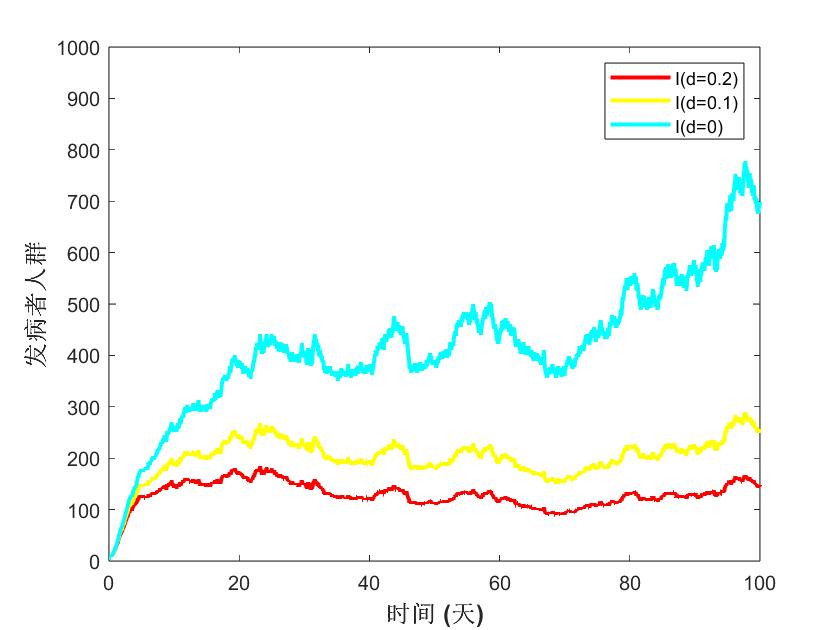} }
	\subfigure[]{
		\includegraphics[width=7.3cm,height=3.9cm]{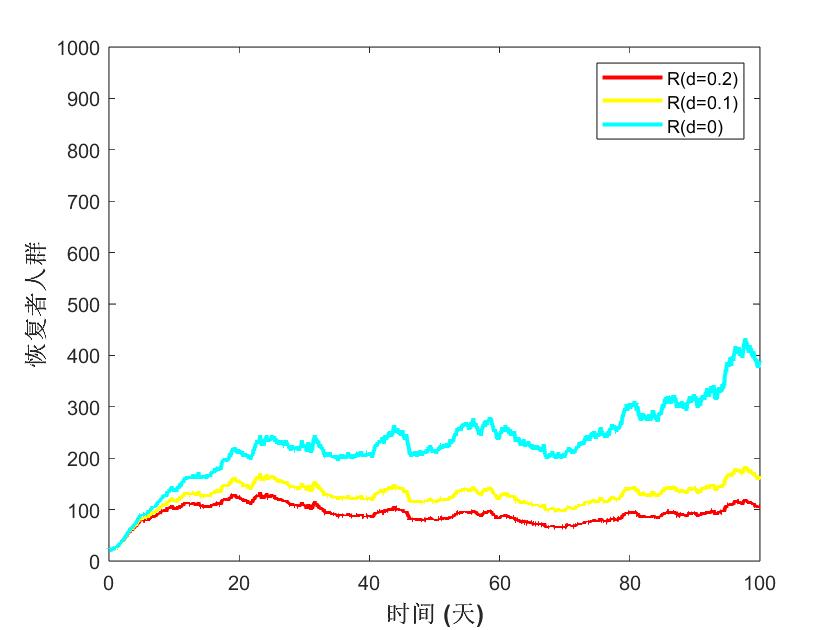} }
\end{minipage}
\vspace{-3mm} 
\caption{ The effect of diseased-death $d$ on the epidemic system (\ref{eq2}).}
\label{B1}
\end{figure}

\end{exam}

\begin{figure}[!htb]	
\centering	
\begin{minipage}[t]{1\textwidth}
	\centering
	\subfigure[]{
		\includegraphics[width=7.3cm,height=3.9cm]{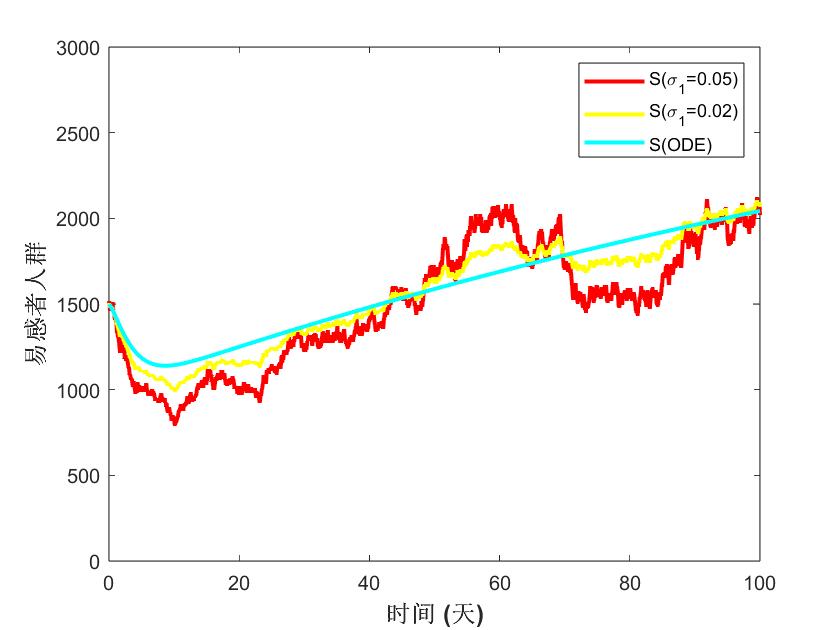} }
	\subfigure[]{
		\includegraphics[width=7.3cm,height=3.9cm]{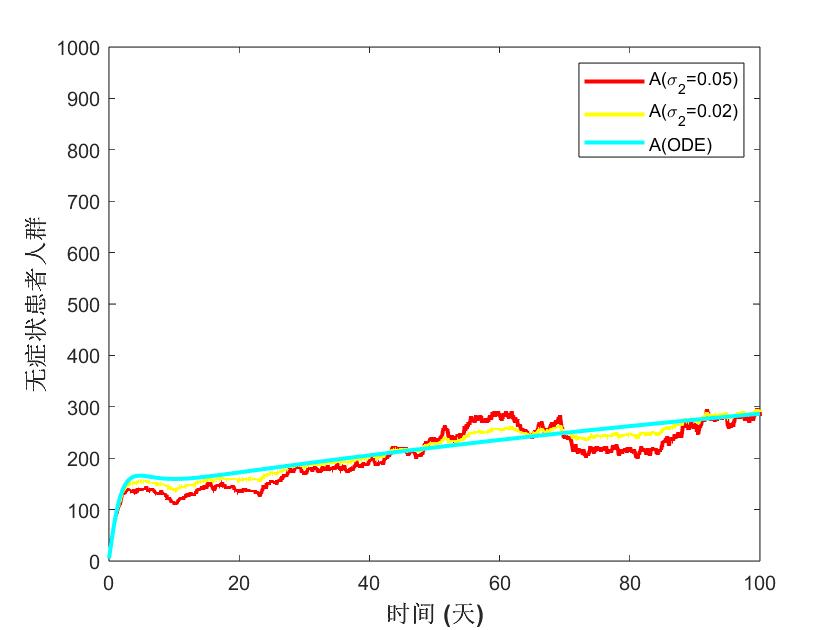} }
\vspace{-4mm} 
\end{minipage}
\begin{minipage}[t]{1\textwidth}
	\centering
	\subfigure[]{
		\includegraphics[width=7.3cm,height=3.9cm]{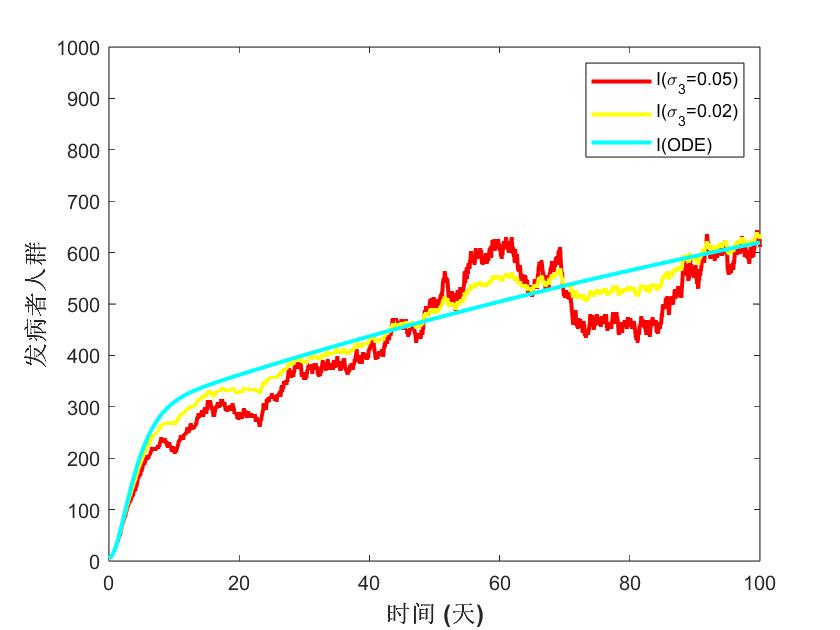} }
	\subfigure[]{
		\includegraphics[width=7.3cm,height=3.9cm]{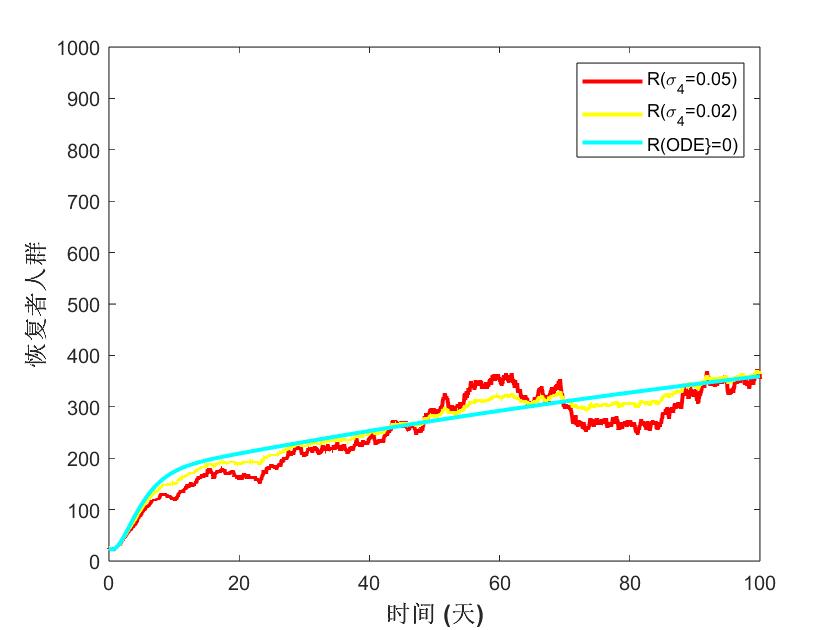} }
\end{minipage}
\vspace{-3mm} 
\caption{ The effect of white noise $\sigma_i(i=1,2,3,4)$ on the epidemic system (\ref{eq2}).}
\label{B2}
\end{figure}

\begin{exam} 
Take $b=0.2$ and $\sigma_i=0, 0.02, 0.05(i=1,2,3,4)$. Other parameter settings are the same as Example \ref{exam1}.  FIGURE \ref{B2} shows the effect of white noise  $\sigma_i$ on the epidemic models (\ref{eq1}) and (\ref{eq2}). Under the influence of white noise, the stochastic epidemic model (\ref{eq2}) fluctuates around the deterministic model (\ref{eq1}), and the magnitude of the fluctuation is related to the degree of $\sigma_{i}(i=1,2,3,4)$.  The size of the fluctuation increases with the value of $\sigma_{i}(i=1,2,3,4)$.

\end{exam}

\begin{exam}
  Take the initial values $S(0)=627000 $,  $A(0)=500 $,  $I(0)=600 $, $R(0)=250000$, and the parameter values $\Lambda=12$, $b=3.1124\times10^{-7} $,  $\beta_{A}=5.0936\times10^{-7} $, $\beta_{I}=5.0725\times10^{-7}$, $\gamma=0.2$, $\mu=2\times10^{-4}$, $\delta_{A}=0.4722$, $\alpha=0.01$, $\delta_{I}=0.9259$, $d=0.0027$, $\sigma_{1}=0.02$, $\sigma_{2}=0.7$, $\sigma_{3}=0.8$, $\sigma_{4}=0.3$, and $\Delta t=0.002 $.  Thus,  $\beta \frac{\Lambda}{\mu}-\frac{h}{2}=-0.6868384<0$, which fulfills the conditions of Theorem \ref{th2.3}.  FIGURE \ref{B3} shows that the disease will eventually disappear according to Theorem \ref{th2.3}.
\begin{figure}[!htb]	
\centering	
\begin{minipage}[t]{1\textwidth}
	\centering
	\subfigure[]{
		\includegraphics[width=7.3cm,height=3.9cm]{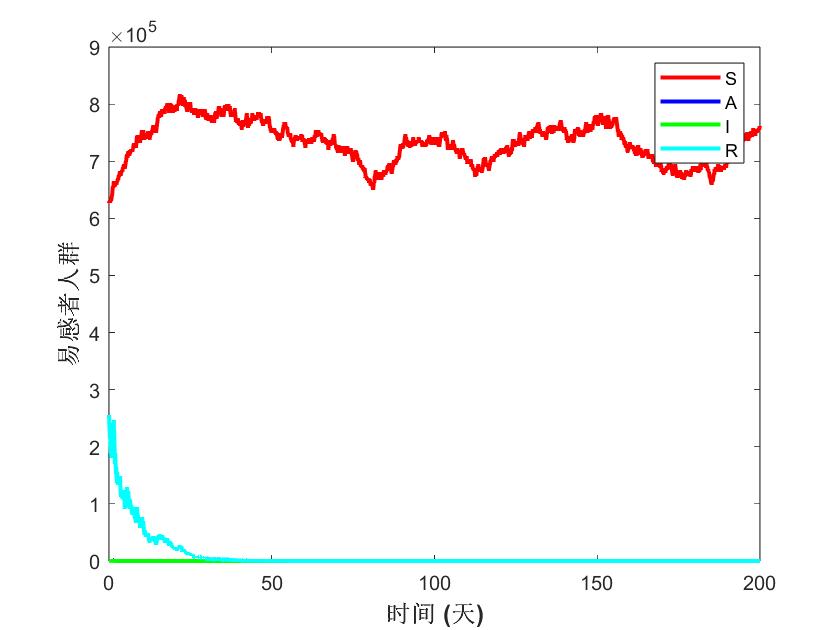} }
	\subfigure[]{
		\includegraphics[width=7.3cm,height=3.9cm]{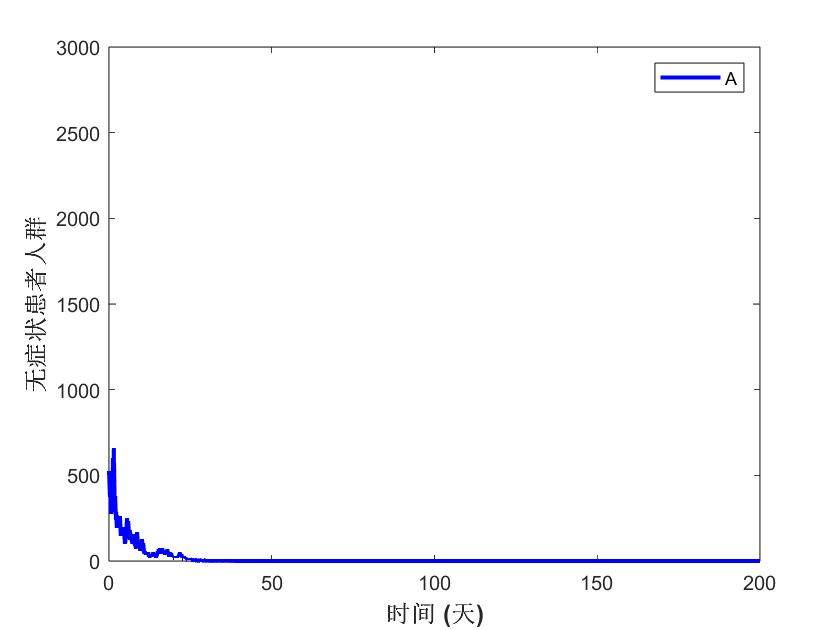} }
 \vspace{-4mm} 
\end{minipage}
\begin{minipage}[t]{1\textwidth}
	\centering
	\subfigure[]{
		\includegraphics[width=7.3cm,height=3.9cm]{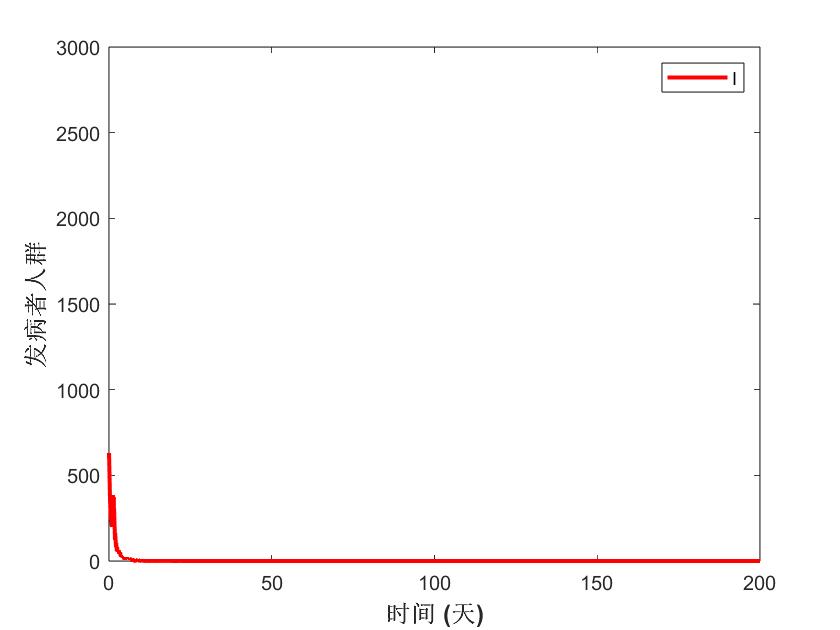} }
	\subfigure[]{
		\includegraphics[width=7.3cm,height=3.9cm]{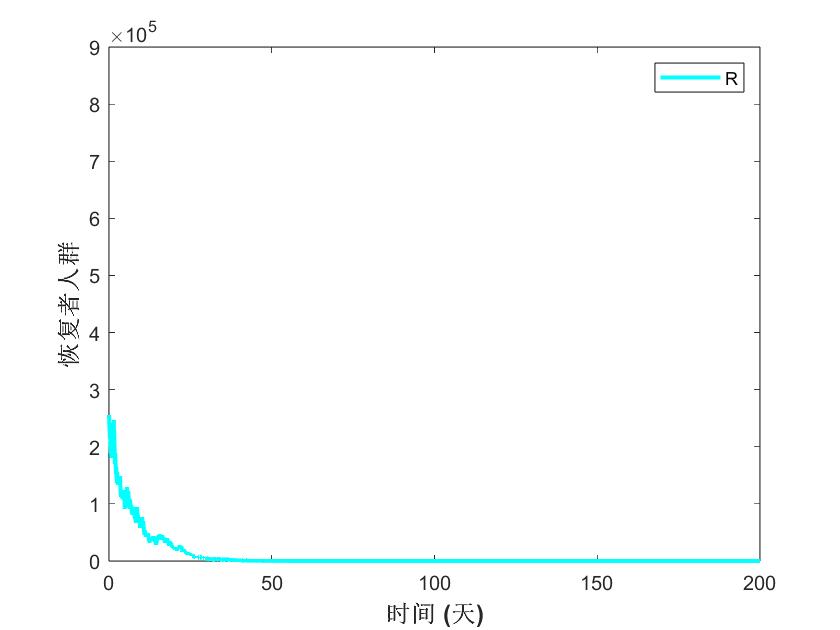} }
\end{minipage}
\vspace{-3mm} 
\caption{ Extinction of the epidemic system  (\ref{eq2}).}
\label{B3}
\end{figure}
\end{exam}
\begin{figure}
\centering	
\begin{minipage}[t]{1\textwidth}
	\centering
	\subfigure[]{
		\includegraphics[width=7.3cm,height=3.9cm]{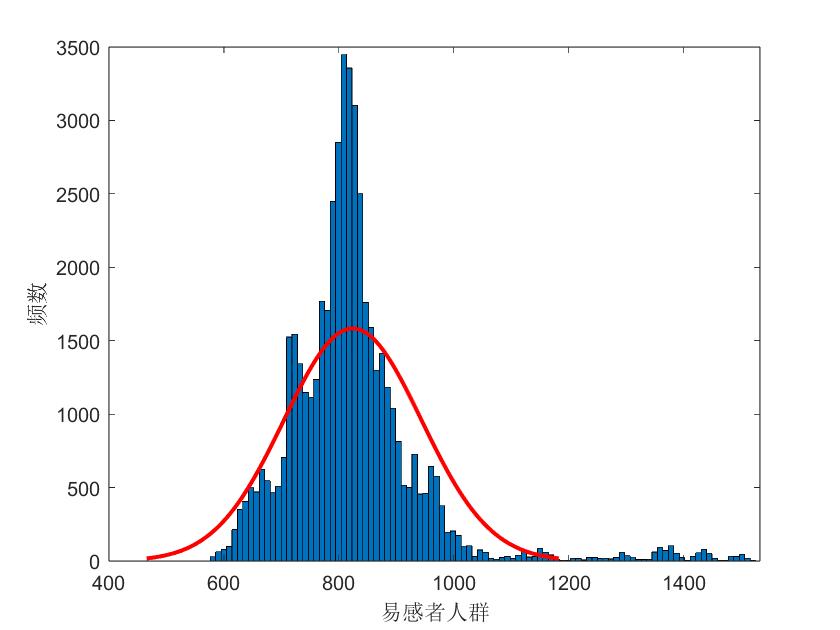} }
	\subfigure[]{
		\includegraphics[width=7.3cm,height=3.9cm]{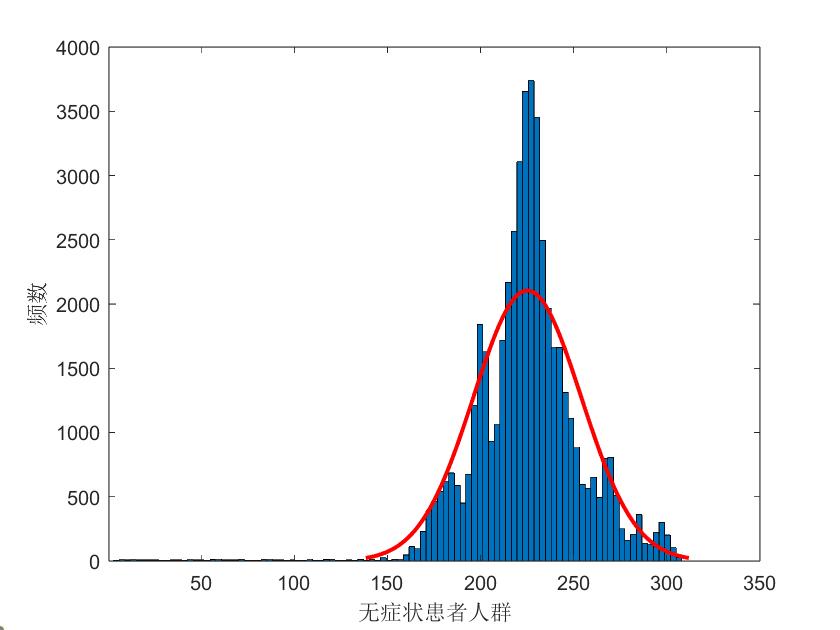} }
 \vspace{-4 mm} 
\end{minipage}
\begin{minipage}[t]{1\textwidth}
	\centering
	\subfigure[]{
		\includegraphics[width=7.3cm,height=3.9cm]{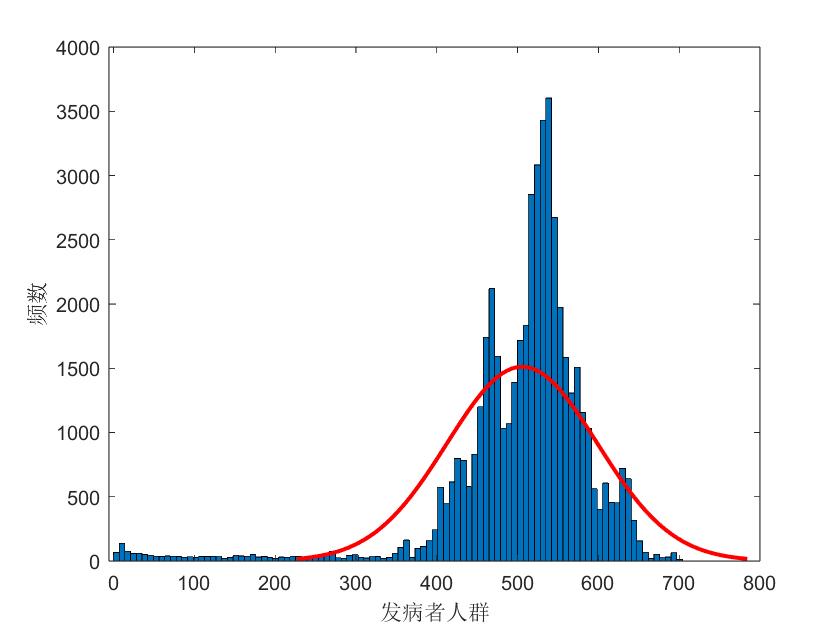} }
	\subfigure[]{
		\includegraphics[width=7.3cm,height=3.9cm]{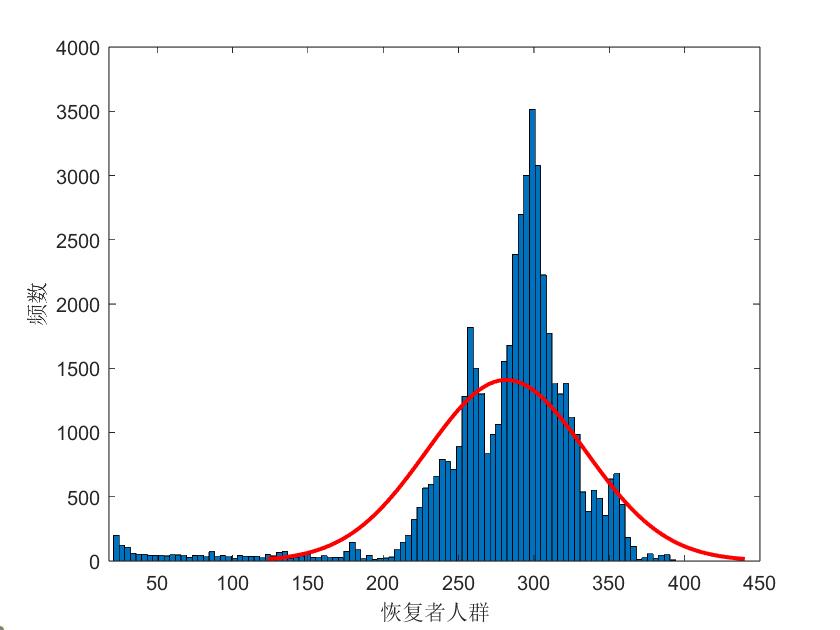} }
\end{minipage}
\vspace{-3mm} 
\caption{ Stationary distribution of system (\ref{eq2}).}
\label{B4}
\end{figure}

\begin{exam}
\textnormal{Take the initial value $S(0)=1500$, $A(0)=5$, $I(0)=6$, $R(0)=25$,  and the parameter values $\Lambda=20 $, $b=0.2 $,  $\beta_{A}=2\times10^{-2} $,  $\beta_{I}=2\times10^{-2} $,   $\gamma=0.5 $,   $\mu=2\times10^{-5} $,  $\delta_{A}=0.2 $,  $\alpha=0.5 $,  $\delta_{I}=0.2$,
 $d=0.0027 $, $\sigma_{1}=0.05 $, $\sigma_{2}=0.05 $,  $\sigma_{3}=0.05 $,  $\sigma_{4}=0.05$, and $\Delta t=0.002 $. Through calcualtion, $R_{0}^{s}=1.16>1$. Based on Theorem \ref{th2.4}, the solution of the model (\ref{eq2}) has a unique stationary distribution. FIGURE \ref{B4} shows the curve and histogram of the stationary distribution for given parameters.}

\end{exam}

\begin{exam}
\textnormal{ Take the initial value $S(0)=624964, A(0)=5, I(0)=6,R(0)=25$, and the parameter values $\Lambda=12, b=4.3\times10^{-5}, \beta_{A}=3\times10^{-5}, \beta_{I}=3\times10^{-5},\gamma=0.5,\mu=7.37\times10^{-3},\delta_{A}=0.92,\alpha=0.072, \delta_{I}=0.92,
d=0.0727,\sigma_{1}=0.1,\sigma_{2}=0.1,\sigma_{3}=0.1,\sigma_{4}=0.1$, and $\Delta t=0.002$. FIGURE \ref{B5} (a) - (d) shows the dynamics of susceptible, asymptomatic, infected, and recovered individuals with and without controls.  From FIGURE \ref{B5} (b) and  (c),  it is evident
that control measures lead to a significant reduction in asymptomatic and infected individuals.  }
\begin{figure}{!htb}
\centering
\begin{minipage}[t]{1\textwidth}
	\centering
	\subfigure[]{
		\includegraphics[width=7.3cm,height=3.9cm]{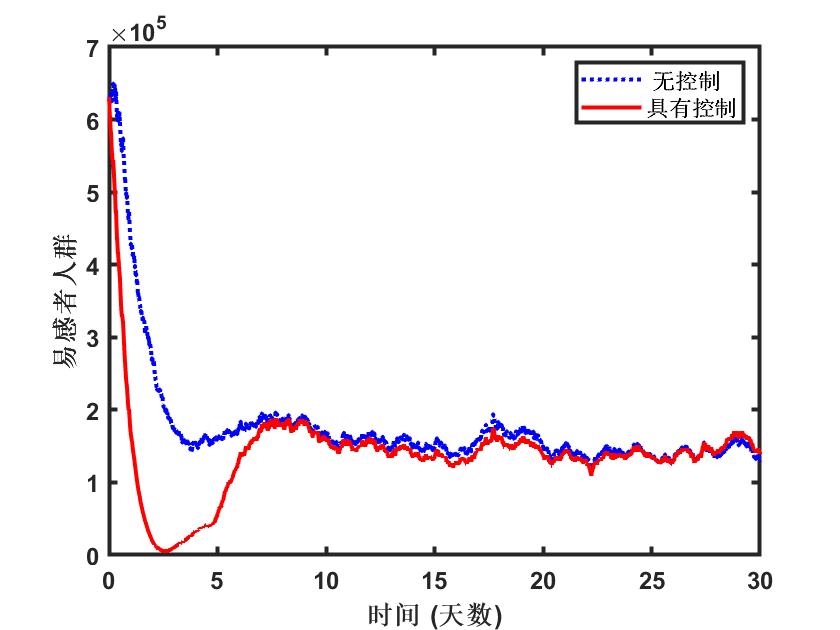} }
	\subfigure[]{
		\includegraphics[width=7.3cm,height=3.9cm]{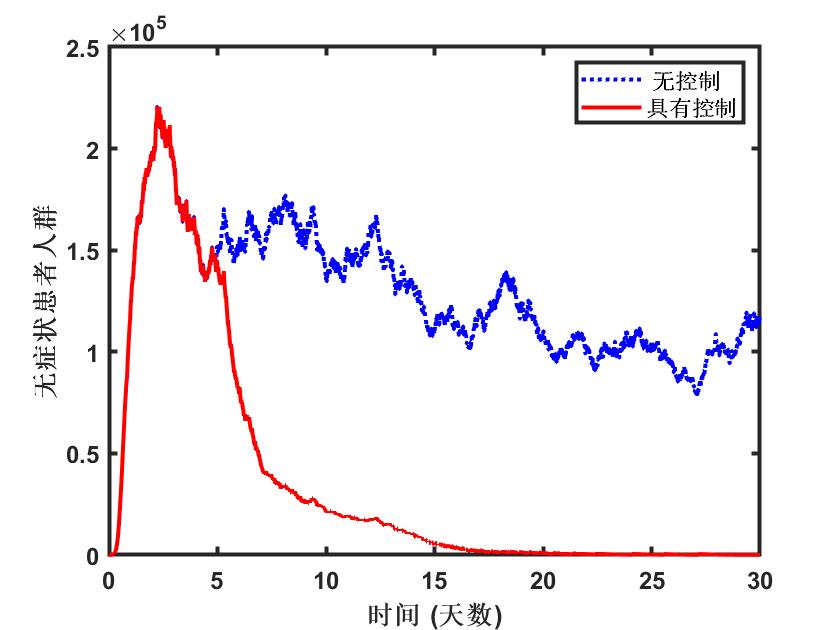} }
\end{minipage}
\begin{minipage}[t]{1\textwidth}
	\centering
	\subfigure[]{
		\includegraphics[width=7.3cm,height=3.9cm]{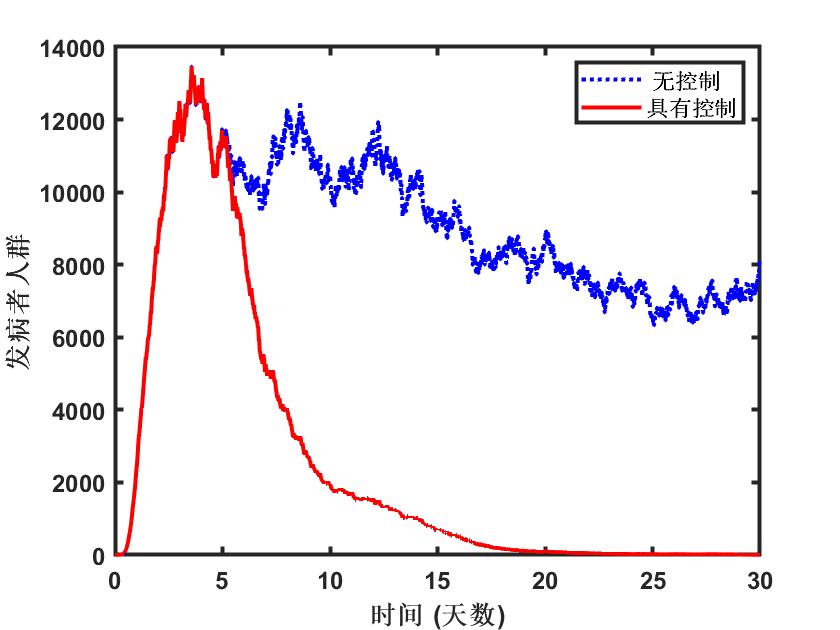} }
	\subfigure[]{
		\includegraphics[width=7.3cm,height=3.9cm]{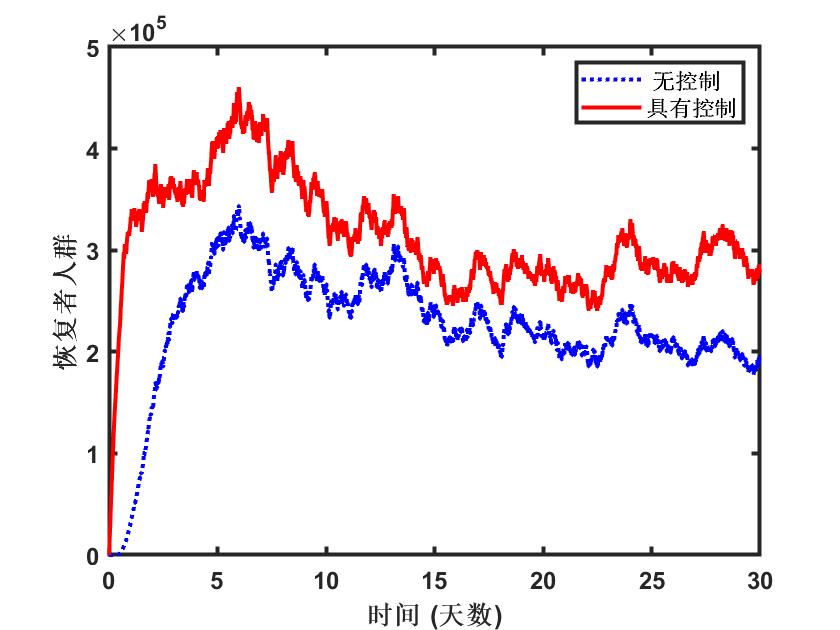} }
\end{minipage}
\vspace{-3mm} 
\caption{ Plots of $S(t),A(t),I(t)$ and $R(t)$ with and without controls in the  system (\ref{eq2}).}
\label{B5}
\end{figure}
\end{exam}

\section{Conclusion}
Based on the deterministic SAIRS model (\ref{eq0}) proposed by \textcolor{blue}{ Ottaviano et al. (2022)}, this paper extends it to a stochastic SAIRS model with saturation incidences and diseased death rate. It proposes a stochastic control system with vaccination and isolation strategies.  Without environmental noise, the stochastic model is equivalent to the deterministic one. That is to say, the model (\ref{eq0}) is a special case of the model (\ref{eq2}) under certain conditions. In the work, we focus on the dynamic properties and optimal control strategies of the model (\ref{eq2}).  Using Lyapunov functions and Ito's formula, the infectious disease is stochastic persistence if $R_{0}^{s}>1$. If $\beta \frac{\Lambda}{\mu}-\frac{h}{2}<0$,  the disease tends to extinction. Further, the system (\ref{eq2}) has a stationary distribution for  $R_{0}^{s}>1$.  Additionally, we introduce two control variables into the model (\ref{eq2}) to determine the optimal control strategy based on Pontryagin's principle of minimization.

Numerical simulation shows the following results: (i) When the saturation coefficient $b$ increases, the susceptible size increases, while asymptomatic, infected, recovered individuals decrease. This means that the behavioral changes and crowing effects of the infected and asymptomatic people can impact the spread of the epidemic. (ii) The diseased death rate $d$ is related to the size of each compartment, especially the infected class. When $d=0,0.1,0.2$, the numbers of $S, A, I, R$ become smaller. It reflects that we need to improve the recovery rate. (iii) The effect of white noises plays a vital role between the models (\ref{eq0}) and (\ref{eq2}). If these noises are small, the dynamic properties of (\ref{eq2}) may be similar to or fluctuate around those of the model (\ref{eq0}). When these noises are significant, it may lead to some different results. (iv) With or without control variables, the model (\ref{eq2}) has different characters. Under an optimal control strategy, the epidemic may be extinct.

 Although we extend the existing model to a stochastic form and study its dynamic properties and control strategies, some problems remain unresolved, including applying or incorporating more complex real-world models.  For example, some scholars have studied the impact of media awareness and global information campaigns \textcolor{blue}{(Das et al. 2020; Taki et al. 2022)}. These open topics are worth further exploration in the future.


\section*{Conflict of interest}
The authors declare no potential conflict of interest.



\nocite{*}
\renewcommand\refname{REFERENCES}

\end{document}